\documentclass[12pt,reqno]{amsart}

\usepackage[arrow,matrix,curve]{xy}

\usepackage[dvips]{graphicx} 

\usepackage{amssymb, latexsym, amsmath, amscd, array,
}

\theoremstyle{definition}

\numberwithin{equation}{section}

\newcommand\N {{\mathbb N}} 
\newcommand\R {{\mathbb R}}
\newcommand\Q {{\mathbb Q}}
\newcommand\Z {{\mathbb Z}} 
\newcommand\st{\text{st}} 

\newcommand\RRR{\mbox{I\!I\!R}}

\newcommand\Los{{\L}o{\'s}}

\newcommand\Horvath{{Horv\'ath}}

\newcommand\fnull{\mathcal{F}_{\!n\!u\!l\!l}}
\newcommand\fez{\mathcal{F}_{ez}}

\DeclareMathOperator{\adequal}{\,{}_{\ulcorner\!\urcorner}\,}

\begin{document}


\thispagestyle{empty}

\title[Leibniz's infinitesimals]{Leibniz's infinitesimals: Their
fictionality, their modern implementations, and their foes from
Berkeley to Russell and beyond}

\author{Mikhail G. Katz and David Sherry}

\address{Department of Mathematics, Bar Ilan University, Ramat Gan
52900 Israel} \email{katzmik@macs.biu.ac.il}

\address{Department of Philosophy, Northern Arizona University,
Flagstaff, AZ 86011, USA} \email{david.sherry@nau.edu}



\bigskip\bigskip\bigskip\noindent
\begin{abstract} 
Many historians of the calculus deny significant continuity between
infinitesimal calculus of the 17th century and 20th century
developments such as Robinson's theory.  Robinson's hyperreals, while
providing a consistent theory of infinitesimals, require the resources
of modern logic; thus many commentators are comfortable denying a
historical continuity.  A notable exception is Robinson himself, whose
identification with the Leibnizian tradition inspired Lakatos,
Laugwitz, and others to consider the history of the infinitesimal in a
more favorable light.  Inspite of his Leibnizian sympathies, Robinson
regards Berkeley's criticisms of the infinitesimal calculus as aptly
demonstrating the inconsistency of reasoning with historical
infinitesimal magnitudes.  We argue that Robinson, among others,
overestimates the force of Berkeley's criticisms, by underestimating
the mathematical and philosophical resources available to Leibniz.
Leibniz's infinitesimals are fictions, not logical fictions, as
Ishiguro proposed, but rather pure fictions, like imaginaries, which
are not eliminable by some syncategorematic paraphrase.  We argue that
Leibniz's defense of infinitesimals is more firmly grounded than
Berkeley's criticism thereof.  We show, moreover, that Leibniz's
system for differential calculus was free of logical fallacies.  Our
argument strengthens the conception of modern infinitesimals as a
development of Leibniz's strategy of relating inassignable to
assignable quantities by means of his transcendental law of
homogeneity.
\end{abstract}

\maketitle

\tableofcontents

\hfill{\small ``Axiom. \; No reasoning about things whereof we have}

\hfill{\small no idea.  Therefore no reasoning about
Infinitesimals.''}

\hfill{\small G. Berkeley, {\em Philosophical commentaries\/}
(no.~354).}

\section{Introduction}

Many historians of the calculus deny any significant continuity
between infinitesimal analysis of the 17th century and non-standard
analysis of the 20th century, e.g., the work of Robinson.  While
Robinson's non-standard analysis constitutes a consistent theory of
infinitesimals, it requires the resources of modern logic;%
\footnote{It is often claimed that the hyperreals require the
resources of model theory.  See Appendix~\ref{rival} for a more
nuanced view.}
thus many commentators are comfortable denying a historical
continuity:
\begin{quote}
[T]here is \dots no evidence that Leibniz anticipated the techniques
(much less the theoretical underpinnings) of modern non-standard
analysis (Earman 1975, \cite[p.~250]{Ea}).
\end{quote}

\begin{quote}
The relevance of current accounts of the infinitesimal to issues in
the seventeenth and eighteenth centuries is rather minimal (Jesseph
1993, \cite[p.~131]{Je93}).
\end{quote}

Robinson himself takes a markedly different position.  Owing to his
formalist attitude toward mathematics, Robinson sees Leibniz as a
kindred soul:
\begin{quote}
Leibniz's approach is akin to Hilbert's original formalism, for
Leibniz, like Hilbert, regarded infinitary entities as ideal, or
fictitious, additions to concrete mathematics.  (Robinson 1967,
\cite[p.~40]{Ro67}).
\end{quote}
The Hilbert connection is similarly reiterated by Jesseph in a recent
text.%
\footnote{\label{f2}See main text at footnote~\ref{f4}.}
Like Leibniz, Robinson denies that infinitary entities are real, yet
he promotes the development of mathematics by means of infinitary
concepts \cite[p.~45]{Ro70}, \cite[p.~282]{Ro66}.  Leibniz's was a
remarkably modern insight that mathematical expressions need not have
a {\em referent\/}, empirical or otherwise, in order to be meaningful.
The fictional nature of infinitesimals was stressed by Leibniz in 1706
in the following terms:
\begin{quote}
Philosophically speaking, I no more admit magnitudes infinitely small
than infinitely great \ldots I take both for mental fictions, as more
convenient ways of speaking, and adapted to calculation, just like
imaginary roots are in algebra. 
%
(Leibniz to Des Bosses, 11 March 1706; in Gerhardt \cite[II,
p.~305]{GP})
\end{quote}

\begin{figure}
\begin{equation*}
\left\{\begin{matrix}assignable \cr quantities\end{matrix}\right\}
\buildrel{\rm LC}\over\leadsto
\left\{\begin{matrix}inassignable \cr quantities\end{matrix}\right\}
\buildrel{\rm TLH}\over\leadsto
\left\{\begin{matrix}assignable \cr quantities\end{matrix}\right\}
\end{equation*}
\caption{\textsf{Leibniz's law of continuity (LC) takes one from
assignable to inassignable quantities, while his transcendental law of
homogeneity (TLH) returns one to assignable quantities.}}
\label{LCLH}
\end{figure}
We shall argue that Leibniz's system for the calculus was free of
contradiction, and incorporated versatile heuristic principles such as
the law of continuity and the transcendental law of homogeneity (see
Figure~\ref{LCLH}) which were, in the fullness of time, amenable to
mathematical implementation as general principles governing the
manipulation of infinitesimal and infinitely large quantities.  And we
shall be particularly concerned to undermine the view that Berkeley's
objections to the infinitesimal calculus were so decisive that an
entirely different approach to infinitesimals was required.

Jesseph suggests an explanation for the irrelevance of modern
infinitesimals to issues in the seventeenth and eighteenth centuries:
\begin{quote}
The mathematicians of the seventeenth and eighteenth centuries who
spoke of taking ``infinitely small" quantities in the course of
solving problems often left the central concept unanalyzed and largely
bereft of theoretical justification (Jesseph \cite[p.~131]{Je93} in
1993).
\end{quote}

According to Jesseph's reading, the earlier approach to infinitesimals
was a conceptual dead-end, and a consistent theory of infinitesimals
required a fresh start.  The force of this claim depends, of course,
on how one understands conceptual analysis and, especially,
theoretical justification.  We argue that Leibniz's defense of the
infinitesimal calculus - both philosophical and mathematical - guided
his successors toward an infinitesimal analysis that is rigorous by
today's standards.  In order to make this case we shall show that
Berkeley's allegations of inconsistency in the calculus stem from
philosophical presuppositions which are neither necessary nor
desirable from Leibniz's perspective.

\section{Preliminary developments}
\label{prior}

A distinction between indivisibles and infinitesimals is useful in
discussing Leibniz, his intellectual successors, and Berkeley.

The term \emph{infinitesimal} was employed by Leibniz in 1673 (see
\cite[series 7, vol. 4, no. 27]{Le2008}).  Some scholars have claimed
that Leibniz was the first to coin the term (e.g., Probst 2008,
\cite[p.~103]{Pr08}).  However, Leibniz himself, in a letter to Wallis
dated 30 march 1699, attributes the term to Mercator:
\begin{quote}
for the calculus it is useful to imagine infinitely small quantities,
or, as Nicolaus Mercator called them, infinitesimals (Leibniz
\cite[p.~63]{Le99}).
\end{quote}
Commentators use the term ``infinitesimal" to refer to a variety of
conceptions of the infinitely small, but the variety is not always
acknowledged.  Boyer, a mathematician and well-known historian of the
calculus writes,
\begin{quote}
In the seventeenth century, however, the infinitesimal and kinematic
methods of Archimedes were made the basis of the differential and the
fluxionary forms of the calculus (Boyer \cite[p.~59]{Boy}).
\end{quote}

This observation is not quite correct.  Archimedes' kinematic method
is arguably the forerunner of Newton's fluxional calculus, but his
infinitesimal methods are less arguably the forerunner of Leibniz's
differential calculus.  Archimedes' infinitesimal method employs {\em
indivisibles\/}.  For example, in his heuristic proof that the area of
a parabolic segment is 4/3 the area of the inscribed triangle with the
same base and vertex, he imagines both figures to consist of
perpendiculars of various heights erected on the base (ibid., 49-50).
The perpendiculars are indivisibles in the sense that they are limits
of division and so one dimension less than the area.  \emph{Qua}
areas, they are not divisible, even if, \emph{qua} lines they are
divisible.  In the same sense, the indivisibles of which a line
consists are points, and the indivisibles of which a solid consists
are planes.  We will discuss the term ``consist of" shortly.

Leibniz's infinitesimals are not indivisibles, for they have the same
dimension as the figures that consist of them.  Thus, he treats curves
as composed of infinitesimal lines rather than indivisible
points. Likewise, the infinitesimal parts of a plane figure are
parallelograms. The strategy of treating infinitesimals as
dimensionally homogeneous with the objects they compose seems to have
originated with Roberval or Torricelli, Cavalieri's student, and to
have been explicitly arithmetized by Wallis (Beeley 2008,
\cite[p.~36ff]{Bee}).

Infinitesimals in this sense occur already in Democritus, who, Edwards
surmises, imagined a cone to consist of infinitesimal triangular
pyramids in order to deduce that its volume is 1/3 base times height
(Edwards 1979, \cite[p.~9]{Ed79}).%
\footnote{Edwards also surmises that Democritus saw that a triangular
pyramid could be completed to form a triangular prism with the same
base and height by adding two more prisms, each with the same base and
height (ibid.).}
Democritus probably used indivisibles too, in arguing that pyramids
with equal heights and bases of equal (but not necessarily congruent)
area have the same volume.%
\footnote{Kepler also used both infinitesimals of the same dimension,
treating a circle, e.g., as consisting of infinitesimal triangles, and
indivisibles, treating an ellipse, e.g., as consisting of its radii
(Boyer 1959, \cite[p.~108-9]{Boy59}).}
In that case Democritus treated a pyramid as consisting of plane
sections, parallel to its base.

Plutarch reports that Democritus raised a puzzle in connection with
treating a cone as consisting of circular sections:
\begin{quote}
If a cone is cut by surfaces parallel to the base, then how are the
sections equal or unequal?  If they were unequal then the cone would
have the shape of a staircase; but if they were equal, then all
sections will be equal, and the cone will look like a cylinder, made
up of equal circles; but this is entirely nonsensical (Plutarch quoted
in Edwards 1979, \cite[p.~8-9]{Ed79}).
\end{quote}
 
This puzzle need not arise for infinitesimals of the same dimension,
with an infinitesimal viewed as a frustum of a cone rather than a
plane section.  Zeno raised a similar but more general puzzle in
connection with treating any continuous magnitude as though it
consists of infinitely many indivisibles.  His metrical paradox
proposes a dilemma: If the indivisibles have no magnitude, then a
figure which consists of them has no magnitude; but if the
indivisibles have some (finite) magnitude, then a figure which
consists of them will be infinite.  Zeno's paradox is, of course, a
puzzle for the idea that a finite magnitude consists of indivisibles.
There is a further puzzle for the idea that a magnitude consists of
indivisibles.  If a magnitude consists of indivisibles, then we ought
to be able to add or concatenate them in order to produce or increase
a magnitude.  But indivisibles are not next to one another; as limits
or boundaries, any pair of indivisibles is separated by what they
limit.  Thus, the concepts of addition or concatenation seem not to
apply to indivisibles.%
\footnote{Difficulties of this sort are what lead thinkers to conceive
of continuous magnitudes kinematically.}

The paradox may not apply to infinitesimals in Leibniz's sense,
however.  For, having neither zero nor finite magnitude, infinitely
many of them may be just what is needed to produce a finite magnitude.
And in any case, the addition or concatenation of infinitesimals (of
the same dimension) is no more difficult to conceive of than adding or
concatenating finite magnitudes.  This is especially important,
because it allows one to represent infinitesimals by means of numbers
and so apply arithmetic operations to them.  This is the fundamental
difference between the infinitary methods of Archimedes (and later
Cavalieri) and the infinitary methods of Leibniz and his followers.

The distinction of indivisible from infinitesimal in Leibniz's sense
is not a difficult one.  Boyer distinguishes proofs by infinitesimal
elements from proofs by indivisibles in various places (e.g., p.~109).
But a failure to keep the distinction before one's mind is a source of
misleading claims about the 17th century calculus.  In what follows,
we shall say that a magnitude consists of infinitesimals just in case
the infinitesimals and the original magnitude have the same dimension.
Otherwise, we shall use the term \emph{indivisible}.

\section{A pair of Leibnizian methodologies}

The existence of separate methodologies in Leibniz was already
apparent to de Morgan, who quipped in 1852 that
\begin{quote}
It is also to be noticed that Leibnitz and the Bernoullis demand the
method of exhaustions, or something equivalent, whenever an objection
is raised to infinitesimals.  They do not face a human enemy with
small shot ; they only use it to kill game (de Morgan 1852,
\cite[p.~324]{de}).
\end{quote}
In his seminal study of Leibniz's methodology, H.~Bos described a pair
of distinct approaches to justifying the calculus:
\begin{quote}
Leibniz considered two different approaches to the foundations of the
calculus; one connected with the classical methods of proof by
``exhaustion", the other in connection with a law of continuity (Bos
1974, \cite[section 4.2, p.~55]{Bos}).
\end{quote}
The first approach relies on an Archimedean ``exhaustion''
methodology.  We will therefore refer to it as the A-methodology.  The
second methodology relies more directly on infinitesimals.  We will
refer to it as the B-methodology, in an allusion to Johann Bernoulli,
who, having learned an infinitesimal methodology from Leibniz, never
wavered from it.%
\footnote{\label{f49} G.~Schubring attributes the first systematic use
of infinitesimals as a foundational concept, to Johann Bernoulli (see
\cite[p.~170, 173, 187]{Sch}).  To note the fact of such systematic
use by Bernoulli is not to say that Bernoulli's foundation is
adequate, or that that it could distinguish between manipulations with
infinitesimals that produce only true results and those manipulations
that can yield false results.  One such infinitesimal distinction
between two types of convergence was provided by Cauchy in 1853 (see
\cite{Ca53}), thereby resolving an ambiguity inherent in his 1821
``sum theorem'' (see Br\aa ting~\cite{Br}; Katz \& Katz \cite{KK11b});
Borovik \& Katz \cite{BK}; B\l aszczyk et al. \cite{BKS}).}

Leibniz's fictionalist attitude toward infinitesimals is no longer
controversial today as it was in the eyes of his closest disciples
such as Bernoulli, l'H\^opital, and Varignon.  The fictional nature of
Leibniz's infinitesimals is clarified by G.~Ferraro in the following
terms:
\begin{quote}
According to Leibniz, imaginary numbers, infinite numbers,
infinitesimals, the powers whose exponents were not ``ordinary"
numbers and other mathematical notions are not mere inventions; they
are auxiliary and ideal quantities that [\ldots] serve to shorten the
path of thought (Ferraro \cite[p.~35]{Fe08};
cf.~Leibniz~\cite[p.~92--93]{Le02}).
\end{quote}
On Ferraro's view, Leibniz's infinitesimals enjoy an {\em ideal\/}
ontological status similar to that of the complex numbers, {\em
surd\/} (irrational) exponents, and other ideal quantities.

Leibniz's 1702 letter to Varignon includes an important enclosure,
recently analyzed by Jesseph \cite{Je11}.  Here Leibniz outlines a
geometrical argument involving quantities~$c$ and~$e$ described as
``not absolutely nothing'', and goes on to comment that~$c$ and~$e$
\begin{quote}
are treated as infinitesimals, exactly as are the elements which our
differential calculus recognizes in the ordinates of curves for
momentary increments and decrements (Leibniz
\cite[p.~104-105]{Le02b}).
\end{quote}
Jesseph argues that Leibniz proposed a pair of methodologies, the
first represented by his {\em Quadratura Arithmetica\/} of 1675, and
the second summarized in the 1702 enclosure.  Jesseph comments that
\begin{quote}
it seems that the infinitesimal here is introduced as something like a
Hilbertian%
\footnote{\label{f4}Jesseph's evocation of Hilbert connects well with
a viewpoint expressed by Robinson (see main text at
footnote~\ref{f2}).}
``ideal element'' that arises when we consider limit cases and seek
what Leibniz termed ``the universality which enables [the calculus] to
include all cases, even that where the given lines disappear''
(Jesseph 2011, \cite[p.~21]{Je11}).
\end{quote}
Jesseph echoes Bos in emphasizing the importance of Leibniz's law of
continuity, described as ``not a mathematical principle, but rather a
general methodological rule with applications in mathematics, physics,
metaphysics, and other sciences'' (Jesseph, ibid.).

Recent Leibniz scholarship branches out into two distinct readings of
Leibnizian infinitesimals.  Bos, Ferraro, \Horvath, Jesseph, and
Laugwitz recognize the presence of a pair of methodologies, namely the
A- and B-methodologies mentioned above.  On the other hand, Arthur,
Levey, and others adopt a {\em syncategorematic\/} reading which
recognizes only the A-methodology, as analyzed in
Subsection~\ref{sync}.  The latter reading, in our view, is due to an
incorrect analysis of Leibniz's fictionalism.

\section{\textbf{\emph{Cum Prodiisset}}}
\label{CP}

Leibniz's text \emph{Cum Prodiisset} \cite{Le01c} (translated by Child
\cite{Ch}) dates from around 1701 according to modern scholars.  The
text is of crucial importance in understanding Leibniz's foundational
stance. We will analyze it in detail in this section.

\subsection{Critique of Nieuwentijt}
\label{41b}

Leibniz begins by criticizing Nieuwentijt, who defended a conception
of infinitesimal according to which the product of two infinitesimals
is always zero.  Mancosu's discussion of Nieuwentijt in \cite[chapter
6]{Ma96} is the only one to date to provide a contextual understanding
of Nieuwentijt's thought.  Leibniz describes Nieuwentijt as

\begin{quote}
being driven to fall back on assumptions that are admitted by no one;
such as that something different is obtained by multiplying 2 by~$m$
and by multiplying~$m$ by 2; that the latter was impossible in any
case in which the former was possible; also that the square or cube of
a quantity is not a quantity or Zero (Leibniz translated by Child
\cite[p. 146]{Ch}).
\end{quote}
Leibniz rejects nilsquare and nilcube infinitesimals,%
\footnote{See further in footnotes \ref{arthur1} and \ref{arthur2}.}
which are altogether incompatible with his approach to differential
calculus, as we will see in Subsection~\ref{assign}.

\subsection{Law of Continuity, with examples}
\label{42}

\emph{Cum prodiisset} erects infinitesimal calculus upon the
foundation of Leibniz's law of continuity (LC).  Because it takes a
variety of forms (cf.~Jorgensen 2009, \cite[p.~224-229]{Jo}), the law
of continuity is perhaps best understood as a \emph{family resemblance
concept} following Wittgenstein~\cite{Wit}, i.e., a cluster of related
concepts.  Here Leibniz formulated LC in the following terms:
\begin{quote}
Proposito quocunque transitu continuo in aliquem terminum desinente,
liceat raciocinationem communem instituere, qua ultimus terminus
comprehendatur (Leibniz \cite[p.~40]{Le01c}).
\end{quote}
The passage can be translated as follows:
\begin{quote}
In any supposed continuous transition, ending in any terminus, it is
permissible to institute a general reasoning, in which the final
terminus may also be included.%
\footnote{Boyer claims that Leibniz used this formulation of LC in ``a
letter to [Pierre] Bayle in 1687'' (Boyer \cite[p.~217]{Boy59}).
Boyer's claim contains two errors.  First, the work in question is not
a letter to Bayle but rather the \emph{Letter of Mr.~Leibniz on a
general principle useful in explaining the laws of nature, etc.}
(Leibniz 1687, \cite{Le87}).  Second, while this letter does deal with
Leibniz' continuity principle, it does not contain the formulation
\emph{In any supposed continuous transition, ending in any terminus,
etc.}; instead, it postulates that an infinitesimal change of input
should result in an infinitesimal change in the output (this principle
was popularized by Cauchy in 1821 as the \emph{definition} of
continuity \cite[p.~34]{Ca21}).  Boyer's erroneous claims have been
reproduced by numerous authors, including M.~Kline
\cite[p.~385]{Kli}.}
\end{quote}
The expression ``final terminus'' refers to the terminus mentioned
earlier which is the ``ending'' of the said transition.  We have
deliberately avoided using the term \emph{limit} in our translation.%
\footnote{This is consistent with Child's translation: ``In any
supposed transition, ending in any terminus, it is permissible to
institute a general reasoning, in which the final terminus may also be
included" \cite[p.~147]{Ch}.  We have reinstated the adjective
\emph{continuous} modifying \emph{transition} (deleted by Child
possibly in an attempt to downplay a perceived logical circularity of
defining LC in terms of continuity itself).  Jorgensen
\cite[p.~228]{Jo} cites Child's translation and claims in footnote~21
that ``this passage says nothing about continuity''.}
Translating \emph{terminus} as \emph{limit} misleadingly suggests the
modern technical meaning of \emph{limit} as a real-valued operation
applied to sequences or functions.%
\footnote{\label{boyer1}One scholar who was so misled was Boyer (see
further in footnote~\ref{boyer2}).}
For perhaps the same reason, Bos comments that
\begin{quote}
the fundamental concepts of the Leibnizian infinitesimal calculus can
best be understood as extrapolations to the actually infinite of
concepts of the calculus of finite sequences.  I use the term
``extrapolation" here to preclude any idea of taking a limit (Bos [18,
p. 13]).%
\footnote{Bos goes on specifically to criticize the Bourbaki's
\emph{limite} wording ``(Leibniz) se tient tr\`es pr\`es du calcul des
diff\'erences, dont son calcul diff\'erentiel se d\'eduit par un
passage \`a la limite" (Bourbaki \cite[p.~208]{Bou}).}
\end{quote}


Subsequent remarks in \emph{Cum Prodiisset} make it clear that
\emph{terminus} encompasses inassignable quantities, for the following
reasons:
\begin{itemize}
\item
If one assumes that \emph{terminus} is assignable, then there is no
justification for applying LC to inassignables, as Leibniz does.%
\footnote{Specifically, Leibniz treats in detail an inassignable
quantity he refers to as \emph{status transitus} (see
Subsection~\ref{status}).}
\item
Reading \emph{terminus} as strictly assignable undermines the
identification of LC as expressed in \emph{Cum Prodiisset}, and LC as
expressed in a february 1702 letter to Verignon, assumed by many
Leibniz scholars (see Subsection~\ref{2feb02}).
\item
A finitistic reading of LC in \emph{Cum Prodiisset} departs from the
the interpretation as given in Knobloch, Laugwitz, Robinson, and
others.
\item
If all entities both in the ``transition'' itself and the
\emph{terminus} are finite, then LC becomes a tautology, inapplicable
to the three examples Leibniz wishes to apply it to (see below),
raising the question why Leibniz would have stated it at all.
\item
Leibniz follows Kepler in exploiting LC (see Kline
\cite[p.~385]{Kli}).  Kepler's famous dictum (originating with
Cusanus) concerning the circle being viewed as an infinitangular
polygon clearly involves infinitary entities.
\item
In a letter to Wallis \cite{Le454}, Leibniz relied on LC so as to
defend his use of the characteristic triangle (i.e., the relation
$ds^2=dx^2+dy^2$ along a curve), similarly involving infinitary
entities.
\end{itemize}

Leibniz gives several examples of the application of his Law of
Continuity.  We will focus on the following three examples.

\begin{enumerate}
\item
In the context of a discussion of parallel lines, he writes:
\begin{quote}
when the straight line BP ultimately becomes parallel to the straight
line VA, even then it converges toward it or makes an angle with it,
only that the angle is then infinitely small \cite[p. 148]{Ch}.
\end{quote}
\item
Invoking the idea that the term equality may refer to equality up to
an infinitesimal error, Leibniz writes:
\begin{quote}
when one straight line is equal to another, it is said to be unequal
to it, but that the difference is infinitely small \cite[p. 148]{Ch}.%
\footnote{Equality up to an infinitesimal is a state of transition
from inequality to equality (this anticipates the law of homogeneity
dealt with in Subsection~\ref{homogeneity}).}
\end{quote}
\item
Finally, a conception of a parabola expressed by means of an ellipse
with an infinitely removed focal point is articulated in the following
terms:
\begin{quote}
a parabola is the ultimate form of an ellipse, in which the second
focus is at an infinite distance from the given focus nearest to the
given vertex \cite[p. 148]{Ch}.
\end{quote}
\end{enumerate}

\subsection{\emph{Souverain principe}}
\label{2feb02}

In a 2~feb.~1702 letter to Varignon, Leibniz formulated the law of
continuity as follows:
\begin{quote}
[\ldots] et il se trouve que les r\`egles du fini r\'eussissent dans
l'infini comme s'il y avait des atomes (c'est \`a dire des
\'el\'ements assignables de la nature) quoiqu'il n'y en ait point la
mati\`ere \'etant actuellement sousdivis\'ee sans fin; et que vice
versa les r\`egles de l'infini r\'eussissent dans le fini, comme s'il
y'avait des infiniment petits m\'etaphysiques, quoiqu'on n'en n'ait
point besoin; et que la division de la mati\`ere ne parvienne jamais
\`a des parcelles infiniment petites: c'est parce que tout se gouverne
par raison, et qu'autrement il n'aurait point de science ni r\`egle,
ce qui ne serait point conforme avec la nature du souverain principe
(Leibniz \cite[p.~93-94]{Le02}).
\end{quote}
Knobloch \cite[p.~67]{Kn02}, Robinson \cite[p.~262]{Ro66},%
\footnote{Robinson's attribution in \cite[p.~262]{Ro66} contains a
misprint: ``(Leibniz [1701])'' should be read as ``(Leibniz
[1702])''.}
Laugwitz \cite[p.~145]{Lau92}, and other scholars identify this
passage as an alternative formulation of the law of continuity, which
can be summarized as follows: \emph{the rules of the finite succeed in
the infinite, and conversely}.%
\footnote{Laugwitz pointed out that this law ``contains an {\em a
priori\/} assumption: our mathematical universe of discourse contains
both finite objects and infinite ones'' (Laugwitz
\cite[p.~145]{Lau92}).  As we have already discussed, identifying
distinct A- and B-methodologies in Leibniz does not require realist
commitments.}

\subsection{\emph{Status transitus}}
\label{status}

We resume our analysis of the law of continuity as formulated in
\emph{Cum Prodiisset}.  Leibniz introduces his next observation by the
clause ``of course it is really true that'', and notes that ``straight
lines which are parallel never meet" \cite[p.~148]{Ch}; that ``things
which are absolutely equal have a difference which is absolutely
nothing" \cite[p.~148]{Ch}; and that ``a parabola is not an ellipse at
all" \cite[p.~149]{Ch}.  How does one, then, account for the examples
of Subsection~\ref{42}?  Leibniz provides an explanation in terms of a
state of transition (\emph{status transitus} in the original Latin
\cite[p. 42]{Le01c}):
\begin{quote}
a state of transition may be imagined, or one of evanescence, in which
indeed there has not yet arisen exact equality \ldots{} or
parallelism, but in which it is passing into such a state, that the
difference is less than any assignable quantity; also that in this
state there will still remain some difference, \ldots{} some angle,
but in each case one that is infinitely small; and the distance of the
point of intersection, or the variable focus, from the fixed focus
will be infinitely great, and the parabola may be included under the
heading of an ellipse \cite[p.~149]{Ch}.
\end{quote}

A state of transition in which ``there has not yet arisen exact
equality" refers to example (2) in Subsection~\ref{42};
``parallelism" refers to example (1); including parabola under the
heading of ellipse is example (3).  

Thus, \emph{status transitus} is subsumed under \emph{terminus},
passing into an assignable entity, but is as yet inassignable.
Translating \emph{terminus} as \emph{limit} amounts to translating it
as an assignable entity, the antonym of the meaning intended by
Leibniz.

The observation that Leibniz's \emph{status transitus} is an
inassignable quantity is confirmed by Leibniz's conceding that its
metaphysical status is ``open to question'':

\begin{quote}
whether such a state of instantaneous transition from inequality to
equality, \ldots{} from convergence [i.e., lines meeting--the authors]
to parallelism, or anything of the sort, can be sustained in a
rigorous or metaphysical sense, or whether infinite extensions
successively greater and greater, or infinitely small ones
successively less and less, are legitimate considerations, is a matter
that I own to be possibly open to question \cite[p.~149]{Ch}.%
\footnote{A syncategorematic expression has no \emph{referential}
function. Thus, the phrase `the present king of France is bald' is a
syncategorematic expression, in that it doesn't refer to any concrete
individual.  A syncategorematic expression serves to reveal logical
relations among those parts of the sentence which are referential.  In
Arthur and Levey's interpretation, the infinitesimal ``serves to
reveal logical relations" by tacitly encoding a quantifier applied to
ordinary real values. But Leibniz clearly does not have real values in
mind when he exploits the term \emph{status transitus}. His
\emph{status transitus} is something between real values of the
variable, on the one hand, and its limiting real value, on the
other. Leibniz's observation that the metaphysical (i.e., ontological)
status of infinitesimals is ``open to question" should apparently have
put to rest any suspicions as to their alleged syncategorematic
nature. After all, if an infinitesimal is merely meant as shorthand
for talking about relations among sets of real values, what is the
point of the lingering doubts expressed by Leibniz as to the
ontological legitimacy of infinitesimals? Certainly the absence of a
concrete individual counterpart of the bald king is a closed and shut
question, rather than being ``open to question". See further in
footnotes \ref{arthur1} and \ref{arthur2}.}
\end{quote}

Yet Leibniz asserts that infinitesimals may be utilized independently
of metaphysical controversies:
\begin{quote}
but for him who would discuss these matters, it is not necessary to
fall back upon metaphysical controversies, such as the composition of
the continuum, or to make geometrical matters depend thereon
\cite[p.~149-150]{Ch}.
\end{quote}
To summarize, Leibniz holds that the inassignable status of
\emph{status transitus} is no obstacle to its effective use in
geometry.  The point is reiterated in the next paragraph:

\begin{quote}
If any one wishes to understand these [i.e. the infinitely great or
the infinitely small--the authors] as the ultimate things, or as truly
infinite, it can be done, and that too without falling back upon a
controversy about the reality of extensions, or of infinite continuums
in general, or of the infinitely small, ay, even though he think that
such things are utterly impossible; it will be sufficient simply to
make use of them as a tool that has advantages for the purpose of the
calculation, just as the algebraists retain imaginary roots with great
profit \cite[p. 150]{Ch}.
\end{quote}
Leibniz has just asserted the possibility of the \emph{mathematical}
infinite: ``it can be done", without \emph{ontological} commitments as
to the reality of infinite and infinitesimal objects.

\subsection{Mathematical implementation of \emph{status transitus}}
\label{45}

We will illustrate Leibniz's concept of \emph{status transitus} by
implementing it mathematically in the three examples mentioned by
Leibniz (see Subsection~\ref{42}).  Example (2) can be illustrated in
terms of a finite positive quantity Leibniz denotes
\[
(d)x
\]
(Bos [18, p. 57] replaced this by~$\underline{\text{d}}x$). The
assignable quantity~$(d)x$ passes via infinitesimal~$dx$ on its way to
absolute~$0$.  Then the infinitesimal~$dx$ is the \emph{terminus}, or
the \emph{status transitus}.  Zero is merely the \emph{shadow} of the
infinitesimal.  This particular \emph{status transitus} is the
foundation rock of the Leibnizian definition of the differential
quotient.

Example (1) of parallel lines can be elaborated as follows. Let's
follow Leibniz in building the line through the point~$(0,1)$ parallel
to the~$x$-axis in the plane.  Line~$L_H$ with~$y$-intercept 1 and
$x$-intercept~$H$ is given by~$y = 1 -\frac{x}{H}$.

Now let~$H$ be infinite. The resulting line~$L_H$ has negative
infinitesimal slope, meets the~$x$-axis at an infinite point, and
forms an infinitesimal angle with the~$x$-axis at the point where they
meet.  We will denote by~$\st(x)$ the assignable (i.e., real) shadow
of a finite~$x$.%
\footnote{The notation ``st" parallels that for the standard part
function in the context of the hyperreals (see Appendix~\ref{rival}).}
Then every finite point~$(x,y)\in L_H$ satisfies
\[
\begin{aligned}
\st(x, y) & = (\st(x), \st(y)) \\& = \left(\st(x), \st\left(1
-\tfrac{x}{H}\right)\right) \\& = (\st(x), 1).
\end{aligned}
\]
Hence the finite portion of~$L_H$ is infinitely close to the
line~$y=1$.  The line~$y=1$ is parallel to the~$x$-axis, and is merely
the \emph{shadow} of the inassignable~$L_H$.  Thus, the parallel line
is constructed by varying the oblique line depending on a
parameter. Such variation comprises the \emph{status transitus}~$L_H$
defined by an infinite value of~$H$.

To implement example (3), let's follow Leibniz in deforming an
ellipse, via a \emph{status transitus}, into a parabola. The ellipse
with vertex (apex) at~$(0,-1)$ and with foci at the origin and
at~$(0,H)$ is given by
\begin{equation}
\label{4.1}
\sqrt{x^2 + y^2} +\sqrt{x^2 + (y - H)^2} = H + 2
\end{equation}
We square \eqref{4.1} to obtain
\begin{equation}
\label{4.2b}
x^2 +y^2 +x^2 +(H-y)^2 + 2\sqrt{(x^2+y^2)(x^2+(H-y)^2)} =H^2+4H+4
\end{equation}
We move the radical to one side
\begin{equation}
\label{4.3}
2\sqrt{(x^2+y^2)(x^2+(H-y)^2)} =H^2+4H+4
-\left(x^2+y^2+x^2+(H-y)^2\right)
\end{equation}
and square again, to obtain after cancellation
\begin{equation}
\label{4.4}
\left(y + 2 + \tfrac{2}{H}\right)^2 - (x^2 + y^2)\left(1 +
\tfrac{4}{H} + \tfrac{4}{H^2}\right) = 0.
\end{equation}
The calculation \eqref{4.1} through \eqref{4.4} depends on the
following habits of \emph{general reasoning} (to echo Child's
translation) with assignable quantities, which are generalized to
apply to inassignable quantities (such as the \emph{terminus/status
transitus}) in accordance with the law of continuity:
\begin{itemize}
\item
squaring undoes a radical;
\item
the binomial formula;
\item
terms in an equation can be transfered to the other side; etc.
\end{itemize}
\emph{General reasoning} of this type is familiar in the realm of
ordinary finite real numbers, but why does it remain valid when
applied to the realm of infinite or infinitesimal numbers? The
validity of transfering such \emph{general reasoning} originally
\emph{instituted} in the finite realm, to the realm of the infinite is
postulated by Leibniz's law of continuity.%
\footnote{When the \emph{general reasoning} being transfered to the
infinite realm is generalized to encompass arbitrary elementary
properties (i.e. first order properties), one obtains the
\Los-Robinson transfer principle (see Appendix~\ref{rival}).}

We therefore apply Leibniz's law of continuity to equation~\eqref{4.4}
for an infinite~$H$.  The resulting entity is still an ellipse of
sorts, to the extent that it satisfies all of the equations
\eqref{4.1} through \eqref{4.4}.  However, this entity is no longer
finite.  It represents a Leibnizian \emph{status transitus} between
ellipse and parabola.  This \emph{status transitus} has foci at the
origin and at an infinitely distant point~$(0,H)$.  Assuming~$x$
and~$y$ are finite, we set~$x_0 = \st(x)$ and~$y_0 = \st(y)$, to
obtain a real shadow of this entity:
\[
\begin{aligned}
\st & \left( \left(y + 2 + \tfrac{2}{H} \right)^2 - (x^2 + y^2)\left(1
+ \tfrac{4}{H} + \tfrac{4}{H^2} \right) \right) = \\& = 
\left(y_0 + 2 + \st \left(\tfrac{2}{H}\right) \right)^2 - \left(x^2_0 +
y^2_0\right) \left(1 + \st \left( \tfrac{4}{H} + \tfrac{4}{H^2}
\right)\right) \\& = (y_0 + 2)^2 - \left(x^2_0 + y^2_0\right) \\& = 0.
\end{aligned}
\]
Simplifying, we obtain
\begin{equation}
\label{43}
y_0= \frac{x_0^2}{4}-1.
\end{equation}
Thus, the finite portion of the \emph{status transitus} \eqref{4.4} is
infinitely close to its \emph{shadow}~\eqref{43}, namely the real
parabola~$y= \frac{x^2}{4}-1$ (in Leibniz's terminology as translated
by Child, ``it is really true'' that this parabola has no focus at
infinity--see Subsection~\ref{status}).  This is the kind of payoff
Leibniz is seeking with his law of continuity.

\subsection{Assignable \emph{versus} unassignable}
\label{assign}

In this section, we will retain the term ``unassignable'' from Child's
translation \cite{Ch} (\emph{inassignabiles} in the original Latin,
see \cite[p.~46]{Le01c}).  After introducing finite quantities~$(d)x,
(d)y, (d)z$, Leibniz notes that
\begin{quote}
the unassignables~$dx$ and~$dy$ may be substituted for them by a
method of supposition even in the case when they are evanescent
\cite[p.~153]{Ch}.
\end{quote}
Leibniz proceeds to derive his multiplicative law in the case~$ay=xv$.
Simplifying the differential quotient, Leibniz obtains
\begin{equation}
\label{46}
\frac{ady}{dx}= \frac{xdv}{dx} + v +dv.
\end{equation}
At this point Leibniz proposes to transfer ``the matter, as we may, to
straight lines that never become evanescent'', obtaining%
\footnote{Child incorrectly transcribes formula~\eqref{47} from
Gerhardt, replacing the equality sign in Gerhardt by a plus sign.
Note that Leibniz himself used the sign~$\adequal$ (see McClenon
\cite[p.~371]{Mc23}).}
\begin{equation}
\label{47}
\frac{a\,(d)y}{dx}= \frac{x\,(d)v}{dx} + v +dv.
\end{equation}
The advantage of \eqref{47} over \eqref{46} is that the expressions
$\frac{(d)y}{dx}$ and~$\frac{(d)v}{dx}$ are assignable (real).
Leibniz points out that ``$dv$ is superfluous''.  The reason given is
that ``it alone can become evanescent''.  The law of homogeneity (see
Subsection~\ref{homogeneity}) is not mentioned explicitly in \emph{Cum
Prodiisset}; therefore the rationale for this step is so far
unsatisfactory.  Discarding the~$dv$ term, one obtains the expected
product formula in this case.  Note that thinking of the left hand
side of~\eqref{47} as the assignable \emph{shadow} of the right hand
side is consistent with Leibniz's example (2) (see
Subsection~\ref{42}).

A final item worth noting is the division by \emph{second}
differentials occurring on page 157:%
\footnote{Child's transcription of formula~\eqref{48} contains
numerous errors: the numerator of the fraction~$\frac{x}{a}$ is
missing; the expression~$\frac{ddv}{ddx}$ appears with a~$y$ in place
of~$v$ in the numerator; the expression~$\frac{2dx\,ddv}{a\,ddx}$
appears with a~$ddx$ in place of~$ddv$ in the numerator.}
\begin{equation}
\label{48}
\frac{ddy}{ddx}= \frac{x}{a} \frac{ddv}{ddx} + \frac{v}{a} +
\frac{2}{a}\frac{dx\,dv}{ddx} + \frac{2dv}{a} +
\frac{2dx\,ddv}{a\,ddx} + \frac{ddv}{a}.
\end{equation}
The final formula on page 158 in Leibniz's text (in Child's
translation) is the assignable version of~\eqref{48}:
\begin{equation}
\label{49}
\frac{ddv}{ddx}= \frac{x}{a} \frac{ddy}{ddx} + \frac{v}{a}+
\frac{2}{a}\frac{dx\,dy}{ddx}.
\end{equation}
Formula~\eqref{49} similarly involves division by second order
differentials.  Division by second order unassignable differentials is
incompatible with the nilsquare approach.%
\footnote{See further in footnotes~\ref{arthur1} and \ref{arthur2}.}

\section{Laws of continuity and homogeneity}
\label{LOC}

As discussed in Section~\ref{CP}, in his 1701 text Leibniz views
parallel lines through the lens of a \emph{terminus}, or \emph{status
transitus}, of intersecting lines forming an infinitesimal angle.

\subsection{From secant to tangent}

A related technique, involved in the determination of the tangent line
from an equation for a secant line, is found in Leibniz's 1684 text
\emph{Nova Methodus}:
\begin{quote}
We have only to keep in mind that to find a \emph{tangent} means to
draw a line that connects two points of the curve at an infinitely
small distance, or the continued side of a polygon with an infinite
number of angles, which for us takes the place of the \emph{curve}
(Leibniz \cite{Le84}; translation from Struik \cite[p.~276]{Stru}).
\end{quote}
Leibniz's final clause here indicates that he viewed a curve as an
infinite-sided polygon with infinitesimal sides.  In the terminology
of Subsection~\ref{status}, the polygon is the \emph{status
transitus}, while the circle itself is merely its assignable
\emph{shadow}.

To elaborate on Leibniz's construction of the tangent, if we have a
formula for a secant line through a pair of variable points~$A,A'$
whose distance~$|AA'|$ tends to zero, the formula remains valid for
the \emph{terminus}, or \emph{status transitus}, when~$|AA'|$ is
infinitesimal (a similar calculation was already detailed in
Subsection~\ref{45}).  Note that the \emph{limits} of~$A$ and~$A'$, in
the modern mathematical sense, are necessarily the \emph{same} point,
making it impossible to build the tangent line.  In other words, the
limit of~$|AA'|$ is~$0$, and a distance of~$0$ doesn't correspond to a
line, but simply to a single point, and is disconnected from the
geometrical notion of tangent.  To understand Leibniz's construction
of the tangent line in \emph{Nova Methodus}, we must steer clear of
\emph{limits} in the modern sense, and rely instead on \emph{terminus}
or \emph{status transitus} as elaborated in \emph{Cum Prodiisset}.%
\footnote{\label{boyer2}S.~L'Huillier (1750--1840) understood
Leibniz's \emph{law of continuity} similarly: ``if a variable quantity
at all stages enjoys a certain property, its limit will enjoy the same
property'' \cite[p.~167]{LH}.  L'Huilier, writing a century before
Weierstrass, is using the term \emph{limit} in its generic sense close
to \emph{terminus}/\emph{status transitus}.  Blinded by the modern
limit doctrine, Boyer comments as follows: ``The falsity of this
doctrine is immediately apparent from the fact that \emph{irrational}
numbers may easily be defined as the limit of sequences of
\emph{rational} numbers, or from the observation that the properties
of a polygon inscribed in a circle are not those of the limiting
figure--the circle'' \cite[p.~256]{Boy59}.  But Boyer's ``limiting
figure'' is an anachronistic imposition, of a post-triumvirate
variety, upon both L'Huilier and Leibniz.  What Leibniz had in mind
was a \emph{status terminus} whose \emph{shadow} is the circle (see
also footnote~\ref{boyer1}).}

\subsection{The arithmetic of infinities in \emph{De Quadratura Arithmetica}}
\label{knob}

Leibniz's ``masterwork on the calculus", {\em De Quadratura
Arithmetica\/}, was written near the end of his stay in Paris
(ca. 1676) (Arthur 2001, \cite[p. 393, note 5]{Le2001}; French
translation in Leibniz \cite{Le2004}).  This text makes it clear that
Leibniz introduced a distinction between equality on the nose, on the
one hand, and approximate equality, or infinite closeness, on the
other.  Knobloch puts it as follows:
\begin{quote}
While up to then two quantities were called equal if their difference
was zero, Leibniz called two quantities equal if their difference can
be made arbitrarily, that is infinitely small \cite[p.~63]{Kn02}.
\end{quote}
The rule governing infinitesimal calculation that Knobloch represents
as Leibniz's rule 2.2, states:
\begin{quote}
2.2~$x, y$ finite,~$x = (y +$ infinitely small)~$\iff$~$x - y \approx
0$ (not assignable difference) (Knobloch \cite[p.~67]{Kn02}).
\end{quote}
Here Knobloch represents the relation of being infinitely close by a
pair of wavy lines.  Concerning Leibniz's rules for the arithmetic of
the infinite, Knobloch comments as follows:
\begin{quote}
In his treatise Leibniz used a dozen rules which constitute his
arithmetic of the infinite.  He just applied them without
demonstrating them, only relying on the {\em law of continuity\/}: The
rules of the finite remains valid in the domain of the infinite
(ibid.).
\end{quote}
Knobloch is alluding to Leibniz's formulation of the law of continuity
in a 2~feb.~1702 letter to Varignon~\cite{Le02} (see
Subsection~\ref{2feb02}).  The last of Leibniz's rules is represented
by Knobloch as rule 12:
\begin{quote}
12:~$x$ divided by~$y$ \quad = \quad ($x$ + infinitely small$_1$)
divided \hbox{$\phantom{mm}$}by ($y$ + infinitely small$_2$) \quad
(Knobloch \cite[p.~68]{Kn02}).
\end{quote}
In other words, rule~12 authorizes a replacement of the right-hand
side, 
\begin{quote}
``($x$+infinitely small$_1$) divided by ($y$ + infinitely
small$_2$)'', 
\end{quote}
by the left-hand side, ``$x$ divided by~$y$'', in infinitesimal
calculations.  As we shall see, Rule 12 is crucial to Leibniz's
conception of the differential quotient,~$dy/dx$.

Thus, to find~$dy/dx$ when~$y=x^2$ one starts with
infinitesimal~$\Delta x$ and forms the infinitesimal quotient
$\frac{\Delta y}{\Delta x}$.  One then simplifies the infinitesimal
quotient, relying on Leibniz's {\em law of continuity\/} (see
Subsection~\ref{42}) to justify each simplification (algebraic
manipulations valid for ordinary numbers, are similarly valid for
infinitesimals), so as to obtain the familiar quantity~$2x+\Delta x$.
Next, to produce the expected answer,~$2x$, for the ``differential
quotient'' (today called the derivative), one applies Leibniz's
Rule~$12$ so as to discard the infinitesimal part~$\Delta x$.%
\footnote{Our analysis of Berkeley's criticism of the proof of the
product rule for differentiation appears in Section~\ref{product}.}
Rule 12 amounts to an application of the transcendental law of
homogeneity (see Subsection~\ref{homogeneity}).

\subsection{Transcendental law of homogeneity}
\label{homogeneity}

Leibniz introduces a law called the {\em transcendental law of
homogeneity\/} (TLH), governing equations involving differentials (as
well as higher-order differentials).  Bos summarizes the law in the
following terms:
\begin{quote}
A quantity which is infinitely small with respect to another quantity
can be neglected if compared with that quantity.  Thus all terms in an
equation except those of the highest order of infinity, or the lowest
order of infinite smallness, can be discarded.  For instance,
\[
a+dx =a
\]
\[
dx+ddy=dx
\]
etc.  The resulting equations satisfy this [\dots] requirement of
homogeneity (Bos \cite[p.~33]{Bos} paraphrasing Leibniz
\cite[p.~381-382]{Le10b}).
\end{quote}
The title of Leibniz's text is \emph{Symbolismus memorabilis calculi
algebraici et infinitesimalis in comparatione potentiarum et
differentiarum, et de lege homogeneorum transcendentali}.  The
inclusion of the transcendental law of homogeneity (\emph{lege
homogeneorum transcendentali}) in the title of the text attests to the
importance Leibniz attached to TLH.

Leibniz's TLH has the effect of eliminating higher-order terms.  The
TLH is also mentioned in Leibniz's \emph{Nova Methodus} \cite{Le84}
(GM V, 224) of 1684 (cf. Bos \cite[p.~33]{Bos}).%
\footnote{Kline opines that ``In response to criticism of his ideas,
Leibniz made various, unsatisfactory replies'' \cite[p.~384]{Kli}, and
proceeds to quote a passage from a letter to Wallis from 30 march 1699
(Kline reports an incorrect year 1690): ``It is useful to consider
quantities infinitely small such that when their ratio is sought, they
may not be considered zero but which are rejected as often as they
occur with quantities incomparably greater [\ldots] Thus if we
have~$x+dx$, $dx$ is rejected [\ldots] Similarly we cannot have $xdx$
and $dx\,dx$ standing together.  Hence if we are to differentiate $xy$
we write $(x+dx)(y+dy)-xy=xdy-xy=xdy+ydx+dx\,dy$.  But here $dx\,dy$
is to be rejected as incomparably less than than~$xdy+ydx$'' (Leibniz
\cite[p.~63]{Le99}).  This summary of the law of homogeneity is
dismissed as ``unsatisfactory'' by Kline.  In fairness it must be
added that Kline wrote two years before the appearance of the seminal
study by Bos \cite{Bos}.}

\subsection{Syncategorematics}
\label{sync}

At variance with Bos' and Jesseph's reading, which accords the
Leibnizian, fictional, infinitesimal the status of a separate approach
as distinct from an Archimedean approach, syncategorematically
inclined scholars maintain that the Leibnizian infinitesimal is merely
shorthand for exhaustion \`a la Archimedes.  This approach originates
with the second, 1990 edition of Ishiguro's book \cite{Is}.  This is a
kind of fictionalism, which Ishiguro describes as \emph{logical
fictionalism}, and we think of as \emph{reductive fictionalism}:
Propositions that refer apparently to fictions may be reduced to
propositions that refer only to standard mathematical entities.  Levey
summarizes the approach as follows:
\begin{quote}
by April of 1676, with his early masterwork on the calculus, {\em De
Quadratura Arithmetica\/}, nearly complete, Leibniz has abandoned any
ontology of actual infinitesimals and adopted the {\em
syncategorematic\/} view of both the infinite and the infinitely small
as a philosophy of mathematics and, correspondingly, he has arrived at
the official view of infinitesimals as {\em fictions\/} in his
calculus (Levey \cite[p.~107]{Le08}).
\end{quote}
According to Levey, Leibnis's fictionalism ``may be styled
Archimedean" \cite[p. 133]{Le08}.  Thus, the syncategorematic reading
seeks to reduce the B-methodology to the A-methodology.  More
precisely, there {\em is\/} no separate B-methodology,
syncategorematically speaking.  Levey argues his claim by citing
Leibniz's comments in the month of april, 1676, even though Leibniz
spent the next forty years publishing mathematics that employed
infinitesimal techniques.%
\footnote{\label{eps}Levey elaborates his position as follows: ``The
syncategorematic analysis of the infinitely small is [\ldots]
fashioned around the order of quantifiers so that only finite
quantities figure as values for the variables.  Thus,
\begin{quote}
(3) the difference~$|a-b|$ is infinitesimal
\end{quote}
does not assert that there is an infinitely small positive value which
measures the difference between~$a$ and~$b$.  Instead it reports,
\begin{quote}
($3^*$) \hbox{For every finite positive value~$\varepsilon$, the
difference~$|a-b|$ is less than~$\varepsilon$.}
\end{quote}
Elaborating this sort of analysis carefully allows one to express the
now-usual epsilon-delta style definitions, etc.''  (Levey
\cite[p.~109-110]{Le08}).  To summarize: no B-methodology,
syncategorematically speaking.  What support does Levey provide for
his nominalistic interpretation?  Leibniz's comments in april, 1676.}

Not content with syncategorematizing infinitesimals right out of
Leibniz's thought, Levey pursues an even more radical thesis
concerning their alleged {\em sudden\/} disappearance:
\begin{quote}
[\dots] within a short few weeks, it's all over for the infinitely
small [\ldots] Good-bye to all the wonderful limit entities: good-bye
parabolic ellipse with one focus at infinity, [\ldots] (Levey
\cite[p.~114-115]{Le08}).
\end{quote}
But did Leibniz indeed bid ``good-bye'' to such ``wonderful
entities''?  In point of fact, Leibniz does refer to just such a
``parabolic ellipse with one focus at infinity'' in 1701, in his text
{\em Cum Prodiisset\/} \cite[p.~46-47]{Le01c}, a quarter of a century
after an exaggerated report of its demise by Levey.  The relevant
passage was cited at the end of Subsection~\ref{42}.

An even more explicit statement of Leibniz's position, in terms of the
\emph{status transitus}, was cited in Subsection~\ref{status}, ending
with the words
\begin{quote} 
the parabola may be included under the heading of an ellipse (Child,
\cite[page 149]{Ch}).
\end{quote}
The idea that ``parabola may be included under the heading of an
ellipse'' is one manifestation of Leibniz's law of continuity, ``a
general methodological rule'' as Jesseph puts it.  The reason Leibniz
gives so much detail here (see the full passage cited in
Subsection~\ref{status}) is readily understood if we assume he is
developing the strategy implicit in a B-methodology.  Namely, to the
extent that he is asserting a non-trivial philosophical or heuristic
principle, he seeks to present a plausible justification concerning
the reliance on fictional objects such as infinitesimals and parabolas
with foci at infinity, and their properties.  On the other hand, if
one assumes syncategorematically that infinitesimals and foci at
infinity are merely shorthand for relations among finite objects
expressed by means of a series of quantifiers, why does Leibniz bother
to formulate, and appeal to, the law of continuity?  If the
syncategorematic interpretation is correct, then the law of continuity
can only be asserting a tautology: a sequence of standard entities
consists of standard entities arranged in a sequence.

Knobloch \cite{Kn02} and Arthur \cite{Ar} claim that Leibniz's
Theorem~6 in \cite{Le1993} (referred to as Leibniz's Proposition~6 in
\cite{Ar}) was a major step toward the development of Riemann sums and
a syncategorematic account of infinitesimals.  Leibniz described the
result as `most thorny' (\emph{spinosissima}).  However, Jesseph
\cite{Je11} points out that Leibniz's argument here involves the
construction of auxiliary curves, and the latter ``requires that we
have a tangent construction that will apply to the original curve''
\cite{Je11}.  The full title of Leibniz's work is \emph{Arithmetical
quadrature of the circle, ellipse and hyperbola}, suggesting that it
is not intended as a general method of quadrature.  In order to
achieve such generality, the transmutation theorem must be invoked,
and that requires a fully general method of tangent construction.
Jesseph notes, however, that
\begin{quote}
although Leibniz's investigations in 1675-76 could show how conic
sections and other well-behaved curves could be handled without
recourse to infinitesimals, he himself understood that there were
limitations to what could be achieved with these methods. Indeed, I
suspect that he set aside the \emph{Arithmetical Quadrature} without
publishing it because he had turned his attention to more powerful
methods that he would introduce in the 1680s in what he called ``our
new calculus of differences and sums, which involves the consideration
of the infinite," and ``extends beyond what the imagination can
attain" (GM 5: 307) (Jesseph~\cite{Je11}).
\end{quote}
The syncategorematic interpretation is at odds with the historical
studies by Bos \cite{Bos} and \Horvath~\cite{Ho86}.  One of the most
salient points is the following.  Leibniz defended his intuitions of
an arbitrary-order infinitesimal against Nieuwentijt's intuition of a
%
%
nilsquare infinitesimal, by passing to the reciprocals and arguing in
support of arbitrary-order infinite quantities.  Whatever the merits of
such an argument, the salient point here is that Leibniz did {\em
not\/} adopt what would have been a more straightfoward and powerful
defense for someone viewing infinitesimals as merely logical fictions,
eliminable by a suitable paraphrase.  Had that been the case, Leibniz
would have pointed out that an infinitesimal is merely a
syncategorematic expression indicative of a logical analysis involving
only finite quantities, in the exhaustive spirit of the A-methodology.
Unless Levey is prepared to declare nilsquare infinitesimals similarly
syncategorematic even in the eyes of Nieuwentijt, Leibniz's response
to the latter furnishes evidence in favor for his commitment to
pursuing a separate B-methodology (see additional details in
Section~\ref{horvath}).

The adjective {\em syncategorematic\/} refers to non-denoting phrases
like `the' and `a', whose logical role was analysed first by medieval
scholastics.  But its application to Leibniz's infinitesimals actually
amounts to a claim that Leibniz anticipated Weierstrass's
$\epsilon,\delta$ techniques.%
\footnote{See footnote~\ref{eps} for Weierstrassian epsilontic details
in Levey.}

Interpreting Leibniz as if he had read Weierstrass already would
appear to fall into the category of feedback-style ahistory criticized
by Grattan-Guinness \cite[p.~176]{Grat04}.  Anticipations of
Weierstrass and Cantor are merely reflections of a philosophical
disposition in favor of a sparse ontology prevalent since the
launching of the great experiment of eliminating infinitesimals from
analysis, favoring a tendentious re-writing of its history.  Exactly
this prompts Mancosu to observe that
\begin{quote}
the literature on infinity is replete with such `Whig' history.
Praise and blame are passed depending on whether or not an author
might have anticipated Cantor and naturally this leads to a completely
anachronistic reading of many of the medieval and later contributions
(Mancosu \cite[p.~626]{Ma09}).
\end{quote}

\subsection{Punctiform and non-punctiform continua}

Critics have objected to historical continuity between Leibniz and
Robinson on the grounds that Leibniz's continuum was non-punctiform,
while Robinson's is punctiform. 

Leibniz does not appear to have thought of the continuum as being made
up of points.  Rather, points merely mark locations on the continuum.
%
%
The view of the continuum as made up of points (therefore ``punctiform
continuum") as later pursued by Cantor, is tied up with a
set-theoretic foundation prevalent in modern mathematics.

On the other hand, Leibniz's mathematics can arguably be
\emph{imbedded} into modern mathematics as was done by Robinson.  The
ingredients introduced, such as the axiom of infinity and the law of
excluded middle, do not appear to contradict what is found in Leibniz
and, on the contrary, provide a basis for a mathematical
implementation of key insights already found in Leibniz, such as the
law of continuity (which became the transfer principle) and the
transcendental law of homogeneity (whose special case became the
standard part function).  The latter seems to have been strangely
ignored by scholars, in spite of the detailed discussion in (Bos 1974,
\cite{Bos}).

Thus, the punctiform nature of the modern approach, be it standard or
non-standard, may be irrelevant to interpreting Leibniz and his
calculus.  Certainly traditional historians acknowledge historical
continuity between the calculus of Newton and Leibniz, on the one
hand, and the calculus of today, on the other, in spite of the
punctiform underpinnings of the latter.%
\footnote{Thus, Boyer writes: ``The traditional view [\ldots] ascribes
the invention of the calculus to [\ldots] Newton and [\ldots]
Leibniz'' \cite[p.~187]{Boy59}.}
This is because the punctiform nature of the modern approach only
comes to the fore in phenomena, not in calculus, but in real analysis
as it emerged in the 19th century.  The same remark applies to
Robinson's approach.

\section{``Marvellous sharpness of Discernment''}
\label{one}

\hfill {\small ``It is curious to observe, what subtilty%
\footnote{Berkeley's spelling.}
and skill this great}

\hfill {\small Genius employs to struggle with an insuperable Difficulty;}

\hfill {\small and through what Labyrinths he endeavours to escape the}

\hfill {\small Doctrine of Infinitesimals''. \quad\quad G. Berkeley,
{\em The Analyst\/}}

\medskip

The title of this section is taken from Berkeley
\cite[Section~XVII]{Be}.  In an analysis of George Berkeley's
criticism as expressed in {\em The Analyst\/}~\cite{Be}, Sherry
\cite{She87} identifies what are actually two distinct criticisms that
are frequently conflated in the literature.%
\footnote{Sherry's dichotomy was picked up by
Jesseph~\cite[p.~124]{Je05}.}
Sherry describes them as
\begin{enumerate}
\item
the metaphysical criticism, and 
\item
the logical criticism
\end{enumerate}
(see \cite[p.~457]{She87}).%
\footnote{Three additional significant aspects of Berkeley's criticism
could be mentioned: (3) a belief in naive indivisibles (this a century
after Cavalieri), i.e., a rejection of infinite divisibility that was
already commonly accepted by mathematicians as early as Wallis and
others; (4) empiricism, i.e., a belief that a theoretical entity is
only meaningful insofar as it has an empirical counterpart, or
referent; this belief ties in with Berkeley's theory of perception
which he identifies with a theory of knowledge; (5) Berkeley's belief
that Newton's attempt to escape a reliance on infinitesimals is futile
(see the epigraph to this Section~\ref{one}).  The latter belief is
contrary to a consensus of modern scholars.  Thus, Pourciau \cite{Pou}
argues that Newton possessed a clear kinetic conception of limit
similar to Cauchy's, and cites Newton's lucid statement to the effect
that ``Those ultimate ratios \dots are not actually ratios of ultimate
quantities, but limits \dots which they can approach so closely that
their difference is less than any given quantity\dots'' See Newton
\cite[p.~39]{New46} and \cite[p.~442]{New99}.  The same point, and the
same passage from Newton, appeared a century earlier in Russell
\cite[item 316, p.~338-339]{Ru03}.}

\subsection{Metaphysical criticism}

Sherry describes Berkeley's {\em metaphysical criticism\/} as
targeting a purportedly contradictory nature of the fluxions and
evanescent increments.  Berkeley similarly targets infinitesimals:
\begin{quote}
The foreign Mathematicians are supposed by some, even of our own, to
proceed in a manner, less accurate perhaps and geometrical, yet more
intelligible. Instead of flowing Quantities and their Fluxions, they
consider the variable finite Quantities, as increasing or diminishing
by the continual Addition or Subduction of infinitely small Quantities
\cite[Section V]{Be}.
\end{quote}
Berkeley's criticism emanates from his philosophical commitment to a
theory of perception,%
\footnote{See Sections~\ref{cluster} and~\ref{varieties} for more detailed
comments on Berkeley's philosophy in relation to infinitesimals.}
anchored in the 18th century empiricist dogma that meaningful
expressions refer to particular \emph{perceptions}.%
\footnote{Note here the similarity between the empiricist dogma and
the syncategorematic interpretation of infinitesimals.  Both assert
that meaningfulness consists ultimately in referring to some favored
type of entity.}
Berkeley's metaphysical criticism is anchored in his empiricist
epistemology, which does not allow infinite divisibility.  This is
illustrated by Berkeley's Question 5, which appears in a long list of
questions at the end of {\em The Analyst\/}:
\begin{quote}
Qu. 5. Whether it doth not suffice, that every assignable number of
Parts may be contained in some assignable Magnitude?  And whether it
be not unnecessary, as well as absurd, to suppose that finite
Extension%
\footnote{\label{extension}Berkeley uses \emph{extension} in the sense
of what we would call today a \emph{continuum}.}
is infinitely divisible?
\end{quote}
This has the unwanted, even absurd consequence that nearly all of
traditional geometry must be abandoned.  It is for this reason that we
claim that Berkeley's criticism of the calculus stands on shakier
grounds than Leibniz's defense.

A number of years later, the metaphysical criticism of infinitesimals
will be expressed most forcefully by Karl Marx.  Marx referred to the
theory of Leibniz and Newton as ``the mystical differential
calculus''.  Echoes of class struggle reverberate through rage and cry
of opposition as Marx notes that Leibniz and Newton
\begin{quote}
believed in the mysterious character of the newly discovered calculus,
that yielded true (and moreover, particularly in the geometrical
application, astonishing) results by a positively false mathematical
procedure. They were thus selfmystified, valued the new discovery all
the higher, enraged the crowd of old orthodox mathematicians all the
more, and thus called forth the cry of opposition, that even in the
lay world has an echo and is necessary in order to pave the way for
something new (Marx \cite[p.~168]{Marx}, cited in Kennedy
\cite[p.~307]{Ken}).
\end{quote}
Struik \cite[p.~187-189]{Stru48} concurred, and Carchedi queried:
``Which view of social reality is hidden behind and informs Marx's
method of differentiation?  Marx differentiates with the eyes of the
social scientist, of the dialectician'' \cite[p.~424]{Car}.  Yet one
can't help wondering whether a, dialectical, elimination of
infinitesimals as a class, with its attendant trading of a simple
algorithmic technique for the manipulation of multiple quantifiers
(and quantifier alernations) genuinely serves the interest of either
the ``lay world'' or the proletariat.

\subsection{Logical criticism}

Meanwhile, Berkeley's {\em logical criticism\/} targets a purported
{\em shift in hypothesis\/}, when a proof
\begin{quote}
is guilty of a {\em fallacia suppositionis\/}, that is, of gaining
certain points by means of one supposition, but subsequently attaining
the final goal by retaining the points just won, in combination with
additional points obtained by replacing the original supposition by
its contradictory (Sherry \cite[p.~257]{She87}).
\end{quote}
Thus, Berkeley's ``ghosts of departed quantities'' criticism of the
definition of the ``differential quotient'' amounted to the following
query: how can a quantity ($dx$) possess a ``ghost'' ($dx\not=0$), and
at the same time be ``departed'' ($dx=0$)?  Alternatively, how can the
infinitesimal analyst have his cake ($dx\not=0$) and eat it, too
($dx=0$)?  Or, as Berkeley colorfully put it,
\begin{quote}
I shall now \ldots observe as to the method of getting rid of such
Quantities, that it is done without the least Ceremony (Berkeley
\cite[Section XVIII]{Be}).
\end{quote}

The {\em exposition\/} of Berkeley's critique of the calculus by
Jesseph~\cite{Je93} is described as ``definitive'' by
Sherry~\cite[p.~127]{She95}, who finds, however, shortcomings in
Jesseph's {\em evaluation\/} of that critique.  Sherry notes
Berkeley's double standard in his atitude toward arithmetic and
geometry.  Thus, Berkeley accepts the practice of arithmetic on the
purely pragmatic grounds of its utility, reflecting an instrumentalist
position.  For Berkeley, arithmetic lacks empirical content, i.e.,
numerals do not denote particular perceptions; the same length can be
3 (feet) or 36 (inches).

Meanwhile, in the case of geometry, including the infinitesimal
calculus, Berkeley adopts a non-instrumentalist approach which insists
upon a subject matter.  Thus, Berkeley criticizes infinitesimals for
possessing no referent, unlike the classical geometry of Euclid, which
does possess a referent, in Berkeley's view.  Sherry notes that
Jesseph
\begin{quote}
doesn't really address the question why Berkeley was reluctant to
explicate the calculus by the tools of his philosophy of arithmetic
and its supporting semantical doctrine that meaningfulness lies in the
use to which terms can be put \cite[p.~127]{She95}.
\end{quote}

Both Berkeley's metaphysical criticism and his logical criticism stem
from philosophical blunders, rather than inherent defects in
infinitesimal calculus.  On the one hand, Berkeley is committed to the
absurd rejection of traditional geometry, and on the other, he
abandons his empiricist dogma as soon as it proves inconvenient for
his purposes.  Berkeley's clericalist agenda may have affected such a
selective application of his empiricist dogma, raising issues of
intellectual integrity already alluded to by de Morgan in his
characteristically caustic style:
\begin{quote}
Dishonesty must never be insinuated of Berkeley. But the Analyst was
intentionally a publication involving the principle of Dr.~Whateley's%
\footnote{The reference is to Richard Whately (1787--1863).}
argument against the existence of Buonaparte; and Berkeley was
strictly to take what he found. The Analyst is a tract which could not
have been written except by a person who knew how to answer it. But it
is singular that Berkeley, though he makes his fictitious character
nearly as clear as afterwards did Whateley, has generally been treated
as a real opponent of fluxions (de Morgan \cite[p.~329]{de}).
\end{quote}

In Section~\ref{product}, we will consider a Berkeleyian response to
the law of homogeneity and the product rule.

\section{Berkeley's 
critique of the product~rule: ``The acme of lucidity''?}
\label{product}

\hfill {\small ``Berkeley attacked the logic of the method of fluxions
or}

\hfill {\small infinitesimal calculus, holding that the infinitesimal
[\ldots]}

\hfill {\small was self-contradictory.  His two ways of bringing this
out}

\hfill {\small are the acme of\; lucidity; \; one concerns\; the
fluxion of a}

\hfill {\small power, the other that of a product.''\; \quad \hfil
J. Wisdom \cite{Wis}}

\medskip

In this section, we analyze Berkeley's critique of the product rule.

\subsection{Berkeley's critique}

Berkeley developed a detailed criticism of the proof of the product
rule for differentiation.  Berkeley refers to
\begin{quote}
{\em Leibnitz\/} and his followers in their {\em calculus
differentialis\/} making no manner of scruple, first to suppose, and
secondly to reject Quantities infinitely small (Berkeley
\cite[Section~XVIII]{Be}) [emphasis in the original].
\end{quote}
Berkeley illustrates the {\em unscrupulous rejection\/}, in the
context of a proof of the product rule, as follows:
\begin{quote}
in the {\em calculus differentialis\/} the main Point is to obtain the
difference of such Product.  Now the Rule for this is got by rejecting
the Product or Rectangle of the Differences.  And in general it is
supposed, that no Quantity is bigger or lesser for the Addition or
Subduction of its Infinitesimal: and that consequently no error can
arise from such rejection of Infinitesimals \cite[Section~XVIII]{Be}.
\end{quote}
Berkeley continues:
\begin{quote}
XIX. And yet it should seem that, whatever errors are admitted in the
Premises, proportional errors ought to be apprehended in the
Conclusion, be they finite or infinitesimal: and that therefore the\,
'\!\!\!${\alpha}\kappa\rho\acute{\iota}\beta\varepsilon\iota\alpha$%
\footnote{Accuracy}
of Geometry requires nothing should be neglected or rejected.
\end{quote}
Here Berkeley is objecting to the last step in the calculation
\begin{equation}
\label{41}
\begin{aligned}
d(uv) &= (u+du)(v+dv)-uv=udv+vdu+du\,dv \\ & =udv+vdu.
\end{aligned}
\end{equation}
Is Berkeley's objection valid?  The last step in
calculation~\eqref{41}, namely
\[
{udv+vdu} + {du\,dv} = {udv+vdu}
\]
is an application of Leibniz's transcendental law of homogeneity%
\footnote{Leibniz had two laws of homogeneity, one for dimension and
the other for the order of infinitesimalness.  Bos states that they
`disappeared from later developments' \cite[p.~35]{Bos}, referring to
Euler and Lagrange.}
(see Bos \cite[p.~33]{Bos}) already mentioned in
Subsection~\ref{knob}.  The law justifies dropping the $du\,dv$ term
on the grounds that, given an equation whose two sides contain
differentials of different orders, one is authorized to discard the
higher-order ones.

In his 1701 text {\em Cum Prodiisset\/} \cite[p.~46-47]{Le01c},
Leibniz presents an alternative justification of the product rule (see
Bos \cite[p.~58]{Bos}).  Here he divides by~$dx$ and argues with
differential quotients rather than differentials.%
\footnote{\label{arthur1}Leibniz freely inverts his infinitesimals,
making it difficult to interpret his infinitesimals in terms of modern
nilsquare ones, as Arthur attempts to do in \cite{Ar} (see also
footnote~\ref{arthur2}).}
Adjusting Leibniz's notation (see Subsection~\ref{assign}) to fit with
the calculation~\eqref{41}, we obtain an equivalent calculation
\[
\begin{aligned}
\frac{d(uv)}{dx} &= \frac{(u+du)(v+dv)-uv}{dx} \\&=
\frac{udv+vdu+du\,dv}{dx} \\&= \frac{udv+vdu}{dx} + \frac{du\,dv}{dx}
\\&= \frac{udv+vdu}{dx}.
\end{aligned}
\]
Under suitable conditions the term~$\frac{du\,dv}{dx}$ is
infinitesimal, and therefore the last step
\begin{equation}
\label{72}
\frac{udv+vdu}{dx} + \frac{du\,dv}{dx} =
u\,\frac{dv}{dx}+v\,\frac{du}{dx}
\end{equation}
is legitimized as an instance of the transcendental law of homogeneity
(see Subsection~\ref{homogeneity}), which authorizes one to discard
the higher-order term in an expression containing infinitesimals of
different orders.

\subsection{Rebuttal of Berkeley's logical criticism}

Berkeley's {\em logical criticism\/} of the calculus is that the
evanescent increment is first assumed to be non-zero to set up an
algebraic expression, and then {\em treated as zero\/} in {\em
discarding\/} the terms that contained that increment when the
increment is said to vanish.

It is open to Leibniz to rebut Berkeley's logical criticism by noting
that the evanescent increment is {\em not\/} ``treated as zero'', but,
rather, merely {\em discarded\/} through an application of the
transcendental law of homogeneity.%
\footnote{This is not to say that Leibniz's system for differential
calculus satisfied modern standards of rigor.  Rather, we are
rejecting the claim by Berkeley and triumvirate historians to the
effect that Leibniz's system contained logical fallacies.}

Jesseph \cite{Je05} quoted approvingly a passage from Cajori who
characterized Berkeley's arguments as ``so many bombs thrown into the
mathematical camp".  On the contrary, Berkeley's criticisms reveal
more about Berkeley's own mathematical and philosophical limitations
than about the shortcomings of the mathematics he attempted to
criticize.  Jesseph discusses a defense of Newton by James Jurin
(1684--1750) in the following terms:
\begin{quote}
On Jurin's analysis, there is no inconsistency in dividing by an
increment {\em o\/} to simplify a ratio and then dismissing any
remaining {\em o\/}-terms as ``vanished" \cite[p.~129]{Je05}.
\end{quote}
Jesseph mentions that Jurin defended Newton ``even to the point of
insisting that a ratio of evanescent increments could subsist even as
the quantifies forming the ratio vanish", but misses the essential
point that, while it is {\em incorrect\/} to say that ``there is no
inconsistency in dividing by an increment {\em o\/} to simplify a
ratio and then dismissing any remaining {\em o\/}-terms as `vanished'
'', if one deletes the last two words from Jurin's phrase, this does
become {\em correct\/}: 
\begin{quote}
there is no inconsistency in dividing by an increment {\em o\/} to
simplify a ratio and then dismissing any remaining {\em o\/}-terms.
\end{quote}

Leibniz could have his cake and eat it, too--but not at Berkeley's
empiricist table, as we discuss in the next section.

\section{``A thorn in their sides''}
\label{cluster}

\hfill {\small ``Hobbes [\ldots] \, sought to mount an empirical
ladder in providing a}

\hfill {\small physicalist explanation of geometry, and had the ladder
pulled out}

\hfill {\small from\, under\, him\, by\, [J.] Wallis.\; \ldots \;
Berkeley,\, in\, his\, criticism\, of}

\hfill {\small Newton, was not so easily routed.''  \quad \quad \quad
E.~Strong \cite[p.~92-93]{strong}}

\medskip 
Strong apparently felt that Berkeley's empiricist ladder was sturdier
than Hobbes'.  Was it?

Berkeley's dual criticism is essentially unanswerable from the
viewpoint of his empiricist logic.  He would have similarly rejected
modern theories of Archimedean continua stemming from the work of
Cantor, Dedekind, and Weierstrass, because they involve infinite
aggregates, as well as rejecting infinitesimal-enriched continua,
since both involve infinitary constructions.  Meanwhile, Berkeley's
empiricism outlaws all infinite objects, as is evident, for instance,
from the following item:
\begin{quote}
Qu. 21. Whether the supposed infinite Divisibility of finite Extension%
\footnote{See footnote~\ref{extension}.}
hath not been a Snare to Mathematicians, and a Thorn in their Sides?
\cite[Question 21]{Be}.
\end{quote}
Once Berkeley's empiricist epistemology is rejected, so are the
obstacles to responding to his pair of criticisms.  A response to the
metaphysical criticism lies in a presentation of a more stratified
(hierarchical) structure, with an A-continuum englobed inside a
B-continuum (see Figure~\ref{31} and Section~\ref{varieties}).

\subsection{Felix Klein on infinitesimal calculus}
\label{klein}

In 1908, Felix Klein described a rivalry of such continua in the
following terms.  Having outlined the developments in real analysis
associated with Weierstrass%
\footnote{Weierstrass's nominalistic reconstruction (as C.~S.~Peirce
called it, and as Burgess \cite{Bu83} might have) was analyzed
in~\cite{KK11a}.}
and his followers, Klein pointed out that
\begin{quote}
The scientific mathematics of today is built upon the series of
developments which we have been outlining.  But an essentially
different conception of infinitesimal calculus has been running
parallel with this [conception] through the centuries
\cite[p.~214]{Kl08}.
\end{quote}
Such a different conception, according to Klein,
\begin{quote}
harks back to old metaphysical speculations concerning the structure
of the continuum according to which this was made up of [...]
infinitely small parts (ibid.).
\end{quote}
Klein appears to imply that those in the Leibnizian tradition are
committed to the reality of infinitesimals.  While this is not true
for Leibniz, at least, Klein's comments do indicate the seriousness
with which {\em some\/} leading mathematicians of the
post-triumvirate%
\footnote{\label{triumvirate}C.~Boyer refers to Cantor, Dedekind, and
Weierstrass as ``the great triumvirate'' (see \cite[p.~298]{Boy}).}
era viewed the tradition of a B-continuum,%
\footnote{To note Klein's appreciation of the tradition of the
B-continuum is not to imply that he was referring to modern theories
thereof; see footnote~\ref{mean} for a fuller discussion of Klein's
views.}
inspite of the ``official'' line as to its alleged banishment by
Cantor, Dedekind, and Weierstrass.

\begin{figure}
\[
\xymatrix@C=95pt{{} \ar@{-}[rr] \ar@{-}@<-0.5pt>[rr]
\ar@{-}@<0.5pt>[rr] & {} \ar@{->>}[d]^{\hbox{st}} & \hbox{\quad
B-continuum} \\ {} \ar@{-}[rr] & {} & \hbox{\quad A-continuum} }
\]
\caption{\textsf{Thick-to-thin: applying the law of homogeneity, or
taking standard part (the thickness of the top line is merely
conventional)}}
\label{31}
\end{figure}

H. Poincar\'e expressed himself similarly in his essay {\em Science
and hypothesis\/}.  Having discussed what is recognizably an
Archimedean continuum, Poincar\'e proceeds to ask the following
question: ``Is the creative power of the mind exhausted by the
creation of the mathematical continuum?"  and concludes: ``No: the
works of Du Bois-Reymond demonstrate it in a striking way".
Poincar\'e then details his position as follows:
\begin{quote}
We know that mathematicians distinguish between infinitesimals of
different orders and that those of the second order are infinitesimal,
not only in an absolute way, but also in relation to those of the
first order.  It is not difficult to imagine infinitesimals of
fractional or even of irrational order, and thus we find again that
scale of the mathematical continuum... (Poincar\'e \cite[p.~50]{Po}).
\end{quote}

\subsection{Did George slay the infinitesimal?}
\label{bell}

\begin{figure}
\includegraphics[height=1.7in]{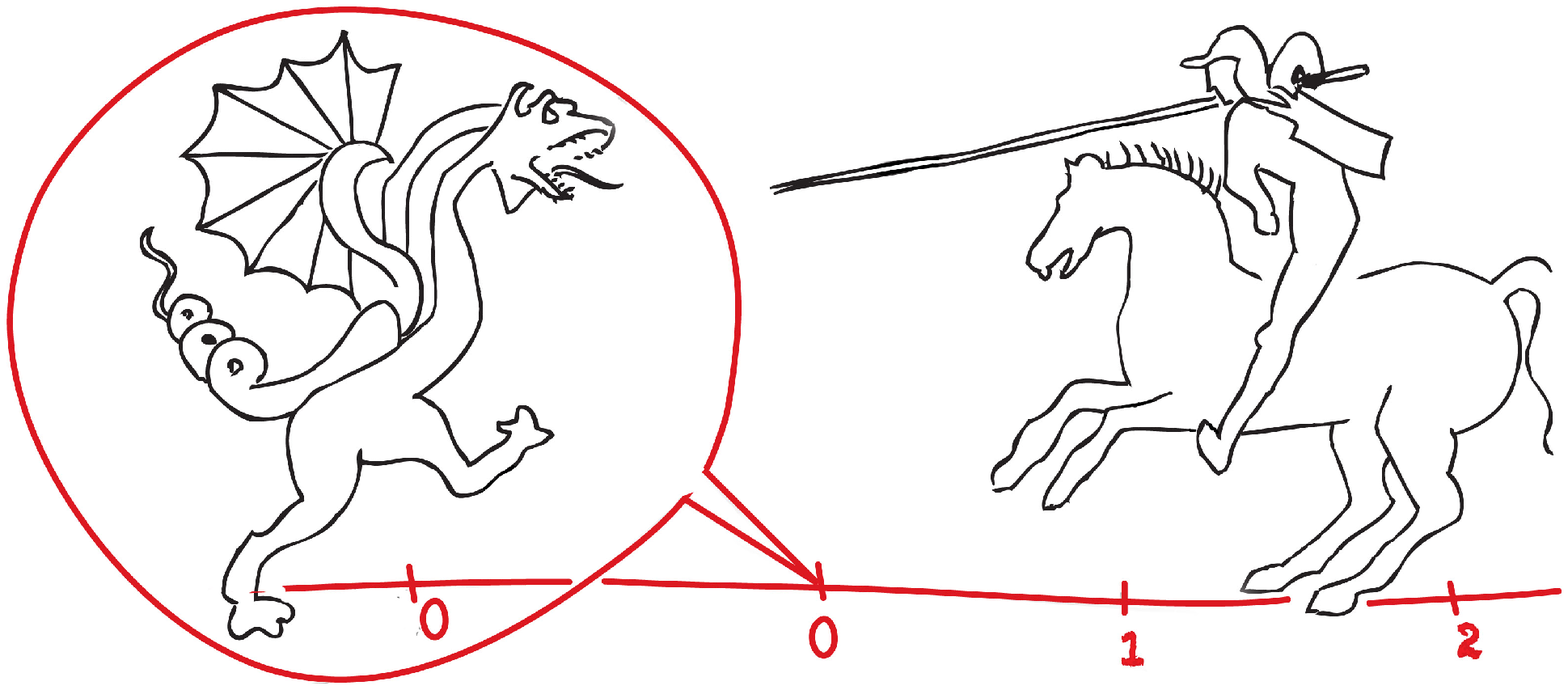}
\caption{\textsf{George's attempted slaying of the infinitesimal,
following E.~T.~Bell and P.~Uccello}}
\label{uccello}
\end{figure}

Not all commentary on the history of the calculus has been as
perceptive as Klein's and Poincar\'e's.  In fact, commentators have
generally regarded Berkeley's criticisms as a crucial step toward
banning infinitesimals from mathematics once and for all; see e.g.,
(Cajori 1917, \cite[p.~151-154]{Caj}) and (Boyer 1949, \cite{Boy}).
In Subsection~\ref{knob} we saw that the means of answering Berkeley's
pair of criticisms were already available to Leibniz himself.  A
historiographic failure to dissect Berkeley's criticism into its two
component parts has led to it being overrated by both historians and
mathematicians, for it prevented them from appreciating the rebuttals
available to Leibniz.  In the course of the 19th century,
infinitesimals came to be thought of as something of an intellectual
embarrassment.  Such a view ultimately found expression in even the
popular (pseudo)historical narratives of E. T. Bell.  Bell waxed
poetic about infinitesimals having been
\begin{itemize}
\item
{\em slain\/} \cite[p.~246]{Bel45},
\item
{\em scalped\/} \cite[p.~247]{Bel45}, and 
\item
{\em disposed of\/} \cite[p.~290]{Bel45} 
\end{itemize}
by the cleric of Cloyne (see Figure~\ref{uccello}).%
\footnote{Dismissing Bell's martial flourishes as merely verbal
excesses would be missing the point.  Bell has certainly been
criticized for other fictional excesses of his purportedly historical
writing (thus, Rothman writes: ``[E. T.] Bell's account [of Galois's
life], by far the most famous, is also the most fictitious''
\cite[p.~103]{Rot}); however, his confident choice of martial imagery
here cannot but reflect Bell's perception of a majority view among
professional mathematicians.  Bell is convinced that Berkeley refuted
infinitesimals only because triumvirate historians and mathematicians
told him so.}
`Scalps' of departed quantities continue to litter the closets of
historical studies of infinitesimal calculus, as we illustrate in
Section~\ref{eleven}.

\section{Varieties of continua}
\label{varieties}

Section~\ref{one} distinguished Berkeley's logical criticism from his
metaphysical criticism.  In 1966 Robinson wrote:
\begin{quote}
The vigorous attack directed by Berkeley against the foundations of
the Calculus in the forms then proposed is, in the first place, a
brilliant exposure of their logical inconsistencies.  But in
criticizing infinitesimals of all kinds, English or continental,
Berkeley also quotes with approval a passage in which Locke rejects
the actual infinite ... It is in fact not surprising that a
philosopher in whose system perception plays the central role, should
have been unwilling to accept infinitary entities
\cite[p.~280-281]{Ro66}.
\end{quote}
Implicit in this passage is Robinson's awareness of the two facets of
Berkeley's criticism.  Unfortunately, he had no access to Leibniz's
manuscript \emph{De Quadratura Arithmetica} first published in 1993
(see Subsection~\ref{knob}).  We therefore cannot agree with his
description of Berkeley's \emph{Analyst} as ``a brilliant exposure''
of ``logical inconsistencies''.  Robinson's praise pays lip service to
the received views on the history of the calculus.  As we showed in
Subsection~\ref{knob}, though, Leibniz possessed the conceptual tools
to formulate a logically unassailable theory of the calculus, in
particular, his use of Rule 12 as a special case of the transcendental
law of homogeneity.

Over the centuries, historians, mathematicians, and philosophers have
envisioned (at least) two distinct theories of the continuum, as
discussed in Section~\ref{cluster}:
\begin{itemize}
\item
an A-continuum (for Archimedes), and
\item
a B-continuum (for Johann Bernoulli,%
\footnote{See footnote~\ref{f49}.}
following Leibniz).
\end{itemize}
The former is a ``thin" continuum, exemplified by what are called
today the real numbers;%
\footnote{\label{f3}Simon Stevin's decimals are at the foundation of
the common number system; see footnote~\ref{stevin} below for
additional details.}
the latter is a ``thick" continuum incorporating infinitesimals.

One possible way of explaning the relation of the two continua is the
following.  All the values in the A-continuum are (theoretically)
possible results of measurement.  The B-continuum has values,
like~$x+dx$, which could never be the result of measurement.

The contents of the A-continuum would correspond, in Leibniz's terms,
to all and only assignable magnitudes.  Meanwhile, the contents of the
B-continuum would contain, in Leibniz's terms, inassignable
magnitudes, as well.%
\footnote{In this context, it may be interesting to note that a close
relationship exists between Cantor's construction of the usual
A-continuum in terms of equivalence classes of Cauchy sequences, on
the one hand, and one of the more straightforward constructions of the
B-continuum, on the other; see Appendix~\ref{rival}.}  The law of
homogeneity as used in~\eqref{72} allows one to pass from a relation
between an inassignable and an assignable quantity, to an equality
between assignable ones.  Note that an equality may be a relation
between two assignable magnitudes, or possibly two inassignable ones.

In the world of L.~Carnot and A.~Cauchy, the assignable/inassignable
distinction translates into a dichotomy of ``constant" quantities
versus ``variable" quantities, an infinitesimal being viewed as
generated by a null sequence%
\footnote{Schubring \cite[p.~454]{Sch} notes that both Cauchy and
Carnot approached infinitesimals dynamically, in terms of sequences
(sometimes referred to as ``variables'', understood as a {\em
succession\/} of values) which tend to zero.  A.~Youschkevitch quotes
Carnot to the effect that an infinitesimal is a variable quantity all
of whose values are determinate and finite (see Gillispie
\cite[p.~242]{Gill}).  To note the fact of the identical definition of
infinitesimals found in Carnot in Cauchy is not to imply total
agreement; thus, Cauchy rejected Carnot's definition of the
differential.}
(see K.~Br\aa ting \cite{Br}).  Cauchy, like most contemporary
authors, did not typically refer to infinitesimals as {\em numbers\/}
because a wide (but not universal) consensus since at least Stevin had
been to use {\em number\/} to refer to the result of counting or
measuring.  The term {\em quantity\/}, as in \emph{quantit\'e
variable}, in Carnot and Cauchy refers to a more general spectrum of
possibilities, typically elements of an {\em ordered\/} system.  Thus,
complex numbers are never referred to as {\em quantities\/} by Cauchy,
but rather as \emph{expressions}.  In Cauchy's presentation of the
material, the starting point are variable quantities, i.e., sequences.
To Cauchy, a sequence ``becomes" an infinitesimal if it tends to zero
(cf.~Borovik \& Katz~\cite{BK}).  More generally, a sequence defines a
variable quantity, that decomposes as a sum of a constant quantity
(i.e. Stevin number) plus a null sequence (that becomes an
infinitesimal).  For example, the ``variable quantity'' defined by the
sequence 
\[
(3.1,\, 3.14,\, 3.141,\, 3.1415,\, \ldots)
\]
will decompose as the sum of the real number $\pi$ and a negative
infinitesimal.  It can easily be shown that a similar decomposition
(for finite elements) necessarily holds in any ordered field properly
containing the reals.

\section{A continuity between Leibnizian and modern infinitesimals?}
\label{seven}

Are modern theories of infinitesimals, legitimate heirs to Leibniz's
theory?  Robinson wrote in 1966:
\begin{quote}
It is shown in this book that Leibniz's ideas can be fully vindicated
and that they lead to a novel and fruitful approach to classical
Analysis and to many other branches of mathematics \cite[p.~2]{Ro66}.
\end{quote}

Such claims on Leibniz's heritage have encountered stiff resistance.
A claim of a historical continuity between Leibniz's and Robinson's
infinitesimals is a formidable challenge to the received triumvirate
scholarship.  There are two ways of deflecting the challenge posed by
modern infinitesimals.  

One of the objections is that Robinson's theory represents a radical
upheaval of the foundations of classical analysis, and the other, than
it is a trivial reformulation thereof.  While these objections are
obviously at odds with each other, they bear scrutiny.  We will
summarize and analyze these objections below.

\begin{enumerate}
\item
There is no historical continuity at all.  Historical infinitesimals
were unsound as analyzed by Berkeley, and became obsolete in 1870.
Robinson's infinitesimals are similar to the historical infinitesimals
in name only.  In reality Robinson's invention relies on sophisticated
model theory and an upheaval of the foundations of mathematics
unimaginable to the pioneers of infinitesimal calculus.
\item
Robinson's approach is merely a re-packaging of the old Weierstrassian
ideas.  If you unwind Robinson's definitions, you find at bottom the
same ideas that revolutionized mathematics in the 1870s, namely the
Cantorian set-theoretic revolution that radically transformed our
understanding of the continuum, from which the modern punctiform
conception has emerged.
\end{enumerate}

While on the face of it, objections (1) and (2) are mutually
contradictory, there is a grain of truth in both, even though both are
off target, as we argue below.

{\em Analysis of criticism (1):} Robinson chose to present his results
in the form of maximal generality, exploiting powerful compactness
theorems (see Malcev \cite{Ma}) to prove the most general existence
results for non-standard extensions of~$\R$.  The price one pays for
generality is that the intuitive source of the notion of an
infinitesimal, namely the null sequence, is well-hidden and buried
deep inside the axiom system, through the introduction of a new
symbol~$\epsilon$ and the system of inequalities
\[
\epsilon < \frac{1}{n}, \quad n\in \N.
\]
On the other hand, the intuitive notion of a null sequence is closer
to the surface in the ultrapower construction (see
Appendix~\ref{rival}).  The latter is less general,%
\footnote{Not every hyperreal field can be obtained this way, though a
more general construction called limit ultrapower can be used to
construct a maximal class hyperreal field (see Ehrlich \cite{Eh12} and
Borovik, Jin, and Katz \cite{BJK}).}
but has the advantage of offering a lucid cognitive and formal link to
historical infinitesimals, generated as they were by variable
quantities or null sequences.

{\em Analysis of criticism (2):} Robinson places himself squarely in
the tradition of classical logic which emerged in the work of Frege,
Peano, and Hilbert at about the same time Cantor, Dedekind, and
Weierstrass banished infinitesimals by reconstructing analysis on the
basis of Stevin numbers.  Thus, Robinson's foundational apparatus is
tame compared to programs proposed by Brouwer and much later Lawvere
\cite{Law}, in the framework of intuitionistic logic.

The crucial point remains that Robinson's use of infinitesimals is in
no way influenced by the fact that they can be defined by means
acceptable to the triumvirate.  Rather, his use of infinitesimals is
governed by Leibniz's law of continuity (see Subsection~\ref{42});
that is, the inferences that Robinson draws by means of infinitesimals
are those that the transfer principle--a precise formulation of the
law of continuity--licenses.

By introducing the hyperreal field, an extension of $\R$, Robinson
gains the problem solving power and convenience of an extended number
system.  Describing such an extension as a trivial modification of
Weierstrassian analysis makes no more sense than claiming that the
latter is a trivial modification of the Greek idea of number, rooted
exclusively in 2, 3, 4, \ldots%
\footnote{From a strict set-theoretic viewpoint, each of the
successive number systems $\N\subset\Z\subset\Q\subset\R$ can be
reduced to the previous one by the familiar set-theoretic
constructions.  Yet it is generally recognized that each successive
enlargement, when additional entities come to be viewed as individuals
(or atomic entities), constitutes a conceptual advantage over the
previous one, with a gain in problem solving power.}
Fields medalist T.~Tao summed up the advantage of the hyperreal
framework by noting that it
\begin{quote}
allows one to rigorously manipulate things such as ``the set of all
small numbers'', or to rigorously say things like ``$\eta_1$ is
smaller than anything that involves $\eta_0$'', while greatly reducing
epsilon management issues by automatically concealing many of the
quantifiers in one's argument (Tao \cite[p.~55]{Tao08}).
\end{quote}

The connection to Weierstrass has to run through each of the,
distinctly un-Weierstrassian, developments listed below.  Namely,
Robinson's theory is anchored in a series of developments each of
which was thought breathtaking in its time:
\begin{enumerate}
\item
(Zermelo 1904, \cite{Ze}) isolated the axiom of choice which until
then was a hidden hypothesis in proofs;
\item
(Tarski 1930, \cite{Tar}) proved the existence of ultrafilters using
the axiom of choice;
\item
(Skolem 1934, \cite{Sk}) constructed non-standard models of Peano
arithmetic;
\item
(Hewitt 1948, \cite{Hew}) constructed hyper-real fields using a form
of the ultrapower construction relying on ultrafilters;
\item
(\Los~1955, \cite{Lo55}) proved his theorem the consequence of which
is the transfer principle for hyper-real fields.
\end{enumerate}

\section{Commentators from Russell onward}
\label{eleven}

\hfill {\small ``There is no such thing as an infinitesimal stretch;
if there were, it}

\hfill {\small would not be an element of the continuum; the Calculus
does not}

\hfill {\small require it,\, and\, to\, suppose\, its\, existence
leads to contradictions.''}

\hfill {\small B. Russell, \emph{The principles of mathematics}, p.~345.}
%

%


\subsection{Russell's {\em non-sequiturs\/}}
\label{russell}

Russell's {\em The principles of mathematics\/} \cite{Ru03} dates from
1903.  Russell opens his discussion of infinitesimals by citing the
Archimedean property:
\begin{quote}
If~$P$,~$Q$ be any two numbers, or any two measurable magnitudes, they
are said to be finite with respect to each other when, if~$P$ be the
lesser, there exists a finite integer~$n$ such that~$nP$ is greater
than~$Q$ \cite[p.~332, start of paragraph 310]{Ru03}.
\end{quote}
This is the first time the word-root of ``measure, measurable" appears
in Russell's chapter XL.  It will play a crucial role in Russell's
argument.  Its meaning needs to be established carefully.  The term
``measurable" can be used in at least the following three senses:
\begin{enumerate}
\item
as a measure-theoretic term, e.g., in the expression {\em
Lebesgue-measurable set\/} (this is the most commonly used sense in
contemporary mathematics);
\item
as a number-theoretic or analytic term, as in ``measurable quantity",
meaning a quantity that can be multiplied by other quantities (such as
the integer $n$ in Russell's comment cited above);
\item
as an empirical term, signifying ``accessible perceptually to the
senses or to physical measuring devices".
\end{enumerate}
The intended meaning here is clearly the number-theoretic meaning~(2).
Russell assumes that $P$ violates the Archimedean property, and notes
that
\begin{quote}
if it were possible for~$Q$ to be \ldots finite,~$P$ would be \ldots
infinitesimal--a case, however, which we shall see reason to regard as
impossible \cite[p.~332, end of paragraph 310]{Ru03}.
\end{quote}
Russell is treading on dangerous ground here, following Cantor's
ill-fated attempt to prove infinitesimals to be inconsistent.  On
page~333, Russell for the first time uses the term ``transfinite" in
place of ``infinite":
\begin{quote}
But it must not be supposed that the ratio of the divisibilities of
two wholes, of which one at least is transfinite, can be measured by
the ratio of the cardinal numbers of their simple parts \cite[p.~333,
paragraph 311]{Ru03}.
\end{quote}
The implication is that the terms ``infinite'' and ``transfinite''
have the same meaning.  In other words, Russell is assuming that any
infinite entity will be a Cantorian infinite entity.  He thus imposes
the straitjacket of Cantor's theory of cardinality upon a discussion
that should should be independent of such theories.%
\footnote{The notion of infinitesimal in Russell's time was a
heuristic concept that has not been defined yet.  The entire
enterprise by Cantor and Russell (to prove the non-well-foundedness of
a heuristic concept that has not been defined yet) retroactively
strikes one as ill-conceived.}

Following a discussion of some \emph{mathematical} uses of the term
``infinitesimal" that Russell finds unobjectionable, Russell makes the
following declaration:
\begin{quote}
What makes these various infinitesimals somewhat unimportant, from a
mathematical standpoint, is, that measurement essentially depends on
the axiom of Archimedes, and cannot, in general, be extended by means
of transfinite numbers \cite[p.~333, paragraph 311]{Ru03}.
\end{quote}

This statement comes at the end of a long paragraph dealing with
purely \emph{mathematical} matters.  What strikes the reader of
Russell's text is the extra-mathematical nature of Russell's
declaration concerning the nature of \emph{measurement}.  Clearly
Russell is now using the root ``measurement" in the empirical
sense~(3).  It is the empirical sense that lends the phrase its
plausibility.  However, as a logical link in Russell's argument, what
is required here is the number-theoretic meaning~(2), as discussed
above.

The conflation of meanings (2) and (3) is not a novelty found in
Russell.  Rather, it can be traced back at least to Berkeley, to whom
no expression is meaningful unless it possesses an empirical
counterpart.

Russell's error is therefore two-fold: first, the identification of
infinity with Cantorian transfinitude, which is understandable to a
certain extent since Cantor's theory was the only theory of infinity
considered reliable by mathematicians at the time; and Russell's
conflation of separate meanings of ``measurable" following Berkeley,
and disregarding several centuries of intervening philosophy, a feat
somewhat less pardonable.

\subsection{Earman and the st-function}
\label{earman}

A few months after Robinson's death J.~Earman published a text
claiming to refute any meaningful connection between Leibniz's and
Robinson's infinitesimals.  Did he succeed?

Earman opens his text by introducing a distinction between two types
of Leibnizian infinitesimals, denoted, respectively, infinitesimal$_1$
and infinitesimal$_2$.

The structural role infinitesimals$_i$,~$i=1,2$, play in Earman's text
is transparent.  Namely, Earman was confronted with overwhelming
evidence that Leibniz thought of infinitesimals as ``ideal" entities.
Earman was seeking to refute such a notion.  Therefore he introduced
the distinction between allegedly two types of infinitesimals:
infinitesimal$_i$,~$i=1,2$.  Such a distinction allowed him to claim
that Leibniz was describing only infinitesimal$_1$ as ``ideal$_j$"
($j=2,3$), but not infinitesimal$_2$.

What is the meaning of Earman's dichotomy?  He says infinitesimal$_1$
is ``intrinsically small", whereas infinitesimal$_2$ ``is incomparably
small with respect to ordinary quantities".%
\footnote{Note that the second sense ties in well with Louis Narens'
approach to measurement where he transforms it into a relative notion:
certain entities are measurable compared to others.  This allows him
to do measurement theory in non-Archimedean contexts (see
\cite{Na76}).}
Is Earman's distinction coherent?  If we do have an absolute/intrinsic
scale of ``ordinary quantities" as Earman puts it, then we also have
an absolute/intrinsic scale of ``smallness", and therefore
infinitesimal$_2$ is englobed in infinitesimal$_1$.

Furthermore, whenever Leibniz deals with an allegedly ``intrinsic"
infinitesimal~$e$, he does not hesitate to consider~$e^2$ and higher
powers, showing that his intrinsically small infinitesimal$_1$ is
englobed in infinitesimal$_2$.  None of the recent Leibniz scholarship
seems to have picked up Earman's infinitesimal$_{i=1,2}$ dichotomy.

On page 239, lines 7-8 Earman claims that Leibniz is referring to
infinitesimal$_1$ in a quote he provides from Leibniz, cited out of
both Loemker \cite{Loe} and Gerhardt \cite{Ge50}.  However, the quote
he provides contains no indication why Leibniz would think these are
allegedly infinitesimals$_1$ rather than infinitesimals$_2$, and
Earman himself does not provide any argument to buttress his claim.

On page 239, lines 15-19 Earman makes a similarly unsubstantiated and
ahistorical claim that Leibniz was not referring to infinitesimals in
an ``ordinal sense", but rather to ``transfinite cardinal numbers".
Such an alleged deduction from Leibniz is a {\em non-sequitur\/}.
Meanwhile, Earman's ahistoral claim plays a key role in his argument.
Thus, on page~240, lines 1-3 he argues that a ``reciprocal of an
infinite cardinal" infinitesimal$_1$ is ``not well defined", and
buttresses his argument in footnote 5 by a reference to Fraenkel's
{\em Abstract Set Theory\/} \cite{Fra} from 1966.  Now such operations
on infinite {\em cardinals\/} are indeed not well defined, but the
assumption that Leibnizian infinitesimals are built from {\em
cardinals\/} rather than from, say, {\em ordinals\/} or other
material, is neither Leibniz's nor Fraenkel's, but rather Earman's.
In point of fact, one can indeed build a consistent theory of
infinitesimals starting with infinite ordinals.  Such a theory is
called the surreals.  P. Ehrlich \cite{Eh12} recently constructed an
isomorphism between maximal surreals and maximal hyperreals.

But Earman's most serious error occurs on page 249, where he discusses
second order infinitesimals in Leibniz, such as~$(dx)^2$.  He
declares:
\begin{quote}
We can always arrange that our non-standard model has second order
infinitesimals (Earman \cite[p.~249]{Ea}).
\end{quote}
Earman goes on to interpret such second-order infinitesimals in a
non-standard model, in his own novel way.  Namely, he seeks to
interpret them as elements of the {\em secondary\/} non-standard
extension, which he denotes by
\[
\R^{**},
\]
of the hyperreal extension~$\R^*$ of~$\R$, so that one
has~$\R\subset\R^*\subset\R^{**}$.  What he means by~$\R^{**}$
is~$(\R^*)^*$.  This can be obtained, for example, by applying the
compactness theorem \cite{Ma}.  Earman's conclusion is that
\begin{quote}
unfortunately, these second-order infinitesimals do not have the
critical property that Leibniz assigned them; namely, if~$\epsilon$ is
a first order infinitesimal, then~$(\epsilon)^2$ is second order
(ibid.).
\end{quote}
Earman is claiming that hyperreal infinitesimals do {\em not\/} have
the property that the square of a first-order infinitesimal is a
second-order infinitesimal.

Now it is true that a square of an infinitesimal in~$\R^*$ will never
yield an element of~$\R^{**}$ which is incomparably smaller than every
element of its subfield~$\R^*\subset\R^{**}$.  But Earman's
dubbing~$\R^{**}$ as being ``second order" is artificial and
unnecessary.  The field~$\R^{**}$ can indeed be described as a
``secondary'' non-standard extension of~$\R$, but this notion of {\em
secondary\/} has nothing to do with second-order infinitesimals!%
\footnote{To elaborate, Earman is simply mistaken to think that
higher-order infinitesimals require higher order hyperreals.  Thus,
given a `first-order' infinitesimal $dx\in\R^*$, a second-order
infinitesimal in Leibniz's sense would correspond to the square
$dx^2\in\R^*$, so that we don't need to consider $\R^{**}$ at all.}
On the other hand, the ordinary square of an infinitesimal in
Robinson's hyperreals~$\R^*$ will indeed be ``second order" in
Leibniz's sense.  Earman's ``second order" criticism of the hyperreals
is merely a play on words and an instance of a strawman fallacy.

Earman grudgingly concedes that
\begin{quote}
It is true that at points non-standard analysis gives the appearance
of following Leibniz's strategy of neglecting infinitesimal terms, but
the appearance is only a superficial one \cite[p.~250]{Ea}.
\end{quote}
What is superficial about such an appearance?  Earman elaborates:
\begin{quote}
In the non-standard definition of derivative, it is true, for example
that the dx terms on the right hand side \ldots is [sic] in a sense
`dropped' in obtaining the answer \ldots But the `dropping' comes from
taking the standard part of the quantities \ldots \cite[p.~250]{Ea}.
\end{quote}
Earman feels that ``dropping'' the~$dx$ terms by Leibniz is dissimilar
from ``dropping'' them by means of the standard part function.  Is it
dissimilar?  Discarding (`dropping') the remaining terms by means of
the standard part function is indeed a modern mathematical
implementation of the transcendental law of homogeneity (see
Subsection~\ref{homogeneity}), following Leibniz's strategy.  Earman
declares that
\begin{quote}
Leibniz's basic strategy of neglecting infinitesimal terms in
comparison with finite ones is not followed in non-standard analysis
\cite[p.~250]{Ea}.
\end{quote}
Earman claims that Leibniz's strategy is not followed in Robinson's
theory.  Or perhaps it is?  Indeed, such ``neglecting" was the object
of Berkeley's logical criticism analyzed in Subsection~\ref{knob} and
Section~\ref{varieties}, where we saw that the binary relation of
equality up to infinitesimal (via the law of homogeneity) was familiar
to Leibniz.  A similar strategy is indeed followed in non-standard
analysis so as to produce a logically consistent definition.  Earman's
remarks amount to criticizing Robinson for providing an answer to
Berkeley's logical criticism.  Earman's remarks are therefore in
error.  He further announces that
\begin{quote}
there is \ldots no evidence that Leibniz anticipated the techniques
\ldots of modern non-standard analysis \cite[p.~250]{Ea}.
\end{quote}
Here Earman contradicts his own discussion on page 245 of Leibniz's
law of continuity, in considering the question of ``which of the
statements true of ordinary quantities are also true of infinitesimal
and infinite quantities, and vice versa?'':
\begin{quote}
Leibniz clearly perceived one form of this question, and in answer he
says that ``the rules of the finite are found to succeed in the
infinite \ldots and conversely, the rules of the infinite apply to the
finite'' (Earman \cite[p.~245]{Ea}).
\end{quote}

Here Earman is citing Leibniz's law of continuity, an antecedent of
the transfer principle (see Subsections~\ref{42}, \ref{2feb02},
and~\ref{knob}).  In this sense it can be said that Leibniz
anticipated a technique of Robinson's theory.

On page~251, in footnote 16, Earman cites Bos' article from 1974, and
writes: ``The reader is urged to consult this excellent article".  Let
us therefore turn to Bos.

\subsection{Bos' appraisal}
\label{bos}

A detailed 1974 scholarly study of Leibniz by Bos contains a brief
appendix dedicated to non-standard issues \cite[Appendix~2]{Bos}.  Bos
notes that
\begin{quote}
to every given function~$f:\R\to\R$, is assigned a unique extension
$f^*: \R^*\to\R^*$, which preserves certain properties of~$f$.  The
field~$\R^*$ provides the framework for the development of the
differential and integral calculus by means of infinitely small and
infinitely large numbers (Bos \cite[p.~81]{Bos}).
\end{quote}
To spell out Bos' suggestion, we choose an infinitesimal
increment~$\Delta x$ and calculate the corresponding increment~$\Delta
y= f^*(x+\Delta x)-f^*(x)$.  Bos \cite[p.~82]{Bos} goes on to
reproduce a version of Robinson's definition of the derivative:
\begin{equation}
\label{101}
f'(x) :=\;^{^{^0}}\!\!\left( \frac{\Delta y}{\Delta x} \right).
\end{equation}
Here the little circle to the upper left of the parenthetical
expression on the right-hand-side of \eqref{101} stands for the
standard part (of the expression included between the parentheses).%
\footnote{See Appendix~\ref{rival}, Subsection~\ref{A2}.  Robinson did
not work with the~``st'' notation, explained in
Sections~\ref{cluster},~\ref{varieties}, and Appendix~\ref{rival}.}
%
%
Bos proceeds to express the following sentiment:%
\footnote{We reproduce Bos' passage in first person.}
\begin{quote}
I do not think that the appraisal of a mathematical theory, such as
Leibniz's calculus, should be influenced by the fact that two and
three quarter centuries later the theory is ``vindicated'' in the
sense that it is shown that the theory can be incorporated in a theory
which is acceptable by present-day mathematical standards (Bos
\cite[p.~82]{Bos}).
\end{quote}
Bos appears to feel that the appraisal of a mathematical theory,
routinely described as logically \emph{inconsistent}, should {\em
not\/} be influenced by a demonstration that it can in fact be
implemented by means of a \emph{consistent} mathematical theory.  Or
perhaps it should be so influenced, as Bos himself pointed out earlier
in his article?  In fact, Bos acknowledged seventy pages earlier that
Robinson's hyperreals do indeed provide~a
\begin{quote}
preliminary explanation of why the calculus could develop on the
insecure foundation of the acceptance of infinitely small and
infinitely large quantities (Bos \cite[p.~13]{Bos}).
\end{quote}
Is {\em that\/} not an instance of an appraisal being influenced by
later vindication?  Furthermore, the clarification of Leibniz's
heuristic law of continuity (see Subsections~\ref{42}
an~\ref{2feb02}): ``rules in the finite domain transfer to the
infinite domain'', in Leibniz's terminology, in terms of a precise
transfer principle provided by \Los's theorem is surely an instance of
a remarkable vindication.  Such a vindication certainly influences our
appreciation of the historical theory.  Thus, Laugwitz points out that
Robinson's infinitesimals
\begin{quote}
can be seen as a revival of Leibniz's {\em fictions\/}, without
explicit mention of their well-foundedness.  The latter is [\ldots]
made precise in a principle of transfer: the {\em rules\/} are founded
on the rules of ``finite totalities'' \cite[p.~152]{Lau92}.
\end{quote}

The most serious error of Bos' examination of Robinson's theory is
that Bos misunderstands the usage of the term ``transfer'' by
Robinson.  Let us analyze the relevant passage from Robinson:
\begin{quote}
Leibniz did say, in one of the passages quoted above, that what
succeeds for the finite numbers succeeds also for the infinite numbers
and vice versa, and this is remarkably close to our {\em transfer of
statements\/} from~$\R$ to~$^*\R$ and in the opposite direction.
(Robinson \cite[p.~266]{Ro66}).
\end{quote}
Bos expresses his disagreement in the following terms:
\begin{quote}
I cannot agree with [Robinson] that this is ``remarkably close to our
transfer of statements from~$\R$ to~$\R^*$ and in the opposite
direction'', and in the context of this passage Robinson himself shows
that Leibniz did not, and could not have, provided such a proof (Bos
\cite[p.~83]{Bos}).
\end{quote}
This may appear to be a disagreement among scholars as to a proper
interpretation of historical theories.  But pay close attention to the
terminology Bos chooses to employ.  Bos quotes Robinson's phrase about
``transfer'', but rejects a suggestion that Leibniz could have
provided ``proof''.  Apparently he understood Robinson's phrase as a
spurious claim of a Leibnizian source for the detailed work of
systematically ``transfering'' (i.e., proving) each of
the~$\R$-statements in an~$\R^*$-context.  None of such detailed work
is to be found in Leibniz, to be sure.  Has Bos then refuted Robinson?
Certainly not.

Indeed, Robinson's reference to ``transfer'' is shorthand for
``transfer principle'' (terminology introduced by later authors),
which makes all real statements that may have been used by Leibniz
{\em automatically\/} true in an~$\R^*$-context.  Robinson was
\emph{not} claiming that Leibniz was busy proving that real statements
also apply to infinitesimals, as it were anticipating such work in
non-standard analysis.  What Robinson did claim is that the heuristic
law of continuity as expressed by Leibniz (see Subsection~\ref{42}),
found expression in a precise metamathematical principle in Robinson's
theory, a point apparently missed by Bos.  The point was not lost on
Urquhart, who commented that
\begin{quote}
some of Bos' criticisms of Robinson involve absurd and impossible
demands - for example, his first criticism (Bos, 1974, p. 83) is that
Robinson proves the existence of his infinitesimals, whereas Leibniz
does not! (Urquhart 2006, \cite{Ur}).
\end{quote}

Bos proceeds to make an apt criticism of Leibniz's theory of
infinitesimals by noting its inability to handle~$\sqrt{h}$ for
infinitesimal~$h$ \cite[p.~83, last line]{Bos}, due to the fact that
Leibniz's infinitesimals come only in integer orders.%
\footnote{Leibniz did consider expressions like~$d^{1/2} x$ and even
gave an explanation of this expression; see his letter to L'Hospital
from 30 september 1695 (Leibniz, \cite{Le95aa}).}
Bos notes that Euler was aware of this problem \cite[p.~84-86]{Bos}.
We may add that Cauchy developed a theory of infinitesimals of
arbitrary (positive) real orders in 1829 (see Laugwitz
\cite[p.~272]{Lau87}), anticipating later theories of non-Archimedean
continua due to Stolz and du Bois-Reymond in terms of rates of growth
of functions, which in turn anticipated Skolem's work on non-standard
models of arithmetic~\cite{Sk} (see the historical discussion in
Robinson \cite[p.~278]{Ro66}).  Some 28 years after Bos, Serfati will
take up the subject of Robinson and the transfer principle again.

\subsection{Serfati's creationist epistemology}
\label{serfati}

In 2002, M.~Serfati challenged Robinson and the transfer principle
in~\cite[p.~317-318]{Se}.  Serfati recognizes that Robinson proved
such a transfer principle {\em between\/} different mathematical
domains, which Serfati describes as a significant accomplishment in
metamathematics.  Serfati proceeds to add, however, that such an
accomplishment does not concern the {\em creation\/} of those
mathematical domains in the first place.  Serfati therefore describes
Robinson's concern as having an ``{\em objectif juridique\/}".
Meanwhile, Serfati's own ``{\em principe de prolongement\/}" is
described as dealing with the {\em creation\/} of mathematical
objects, and therefore is touted as being concerned with the
epistemology of (mathematical) creation, as opposed to what he
describes as ``{\em du juridique\/}", i.e., Robinson's contribution.

Serfati feels that the novelty in this particular area is the {\em
creation\/} of the mysterious entities such as infinitesimals and
infinite numbers.  Once they are created, we can worry about
``juridical'' issues of comparing them with the usual entities, issues
of secondary importance according to Serfati.

The essential point Serfati misses is that non-Archimedean
mathematical objects themselves were already old hat by the time
Robinson came on the scene.  A long line of work on non-Archimedean
systems (Stolz, du Bois-Reymond, Levi-Civita, Hilbert, Skolem, Hewitt,
Schmieden--Laugwitz, to name only a few) testifies to Serfati's
misconception.  The real novelty resides in the transfer principle
allowing one productively to apply such systems in mathematics.%
\footnote{\label{mean}It is interesting to note a criterion of success
of a theory of infinitesimals as proposed by Adolf Abraham Fraenkel
and, before him, by Felix Klein.  In 1908, Klein formulated a
criterion of what it would take for a theory of infinitesimals to be
successful.  Namely, one must be able to prove a mean value theorem
for arbitrary intervals, including infinitesimal ones
\cite[p.~219]{Kl08}.  In 1928, A.~Fraenkel \cite[pp.~116-117]{Fran}
formulated a similar requirement in terms of the mean value theorem.
Such a Klein-Fraenkel criterion is satisfied by the
Hewitt-\Los-Robinson theory by the transfer principle (see
Appendix~\ref{rival}).}

On the last line on page 317, Serfati mentions Leibniz's principle of
continuity, implying that it fits better with his ``creationist
epistemology'' than with Robinson's ``juridical metamathematics''.  In
fact, the opposite is true: Leibniz's law of continuity (see
Subsection~\ref{42}) was precisely a principle relating two domains:
the finite and the infinite.  Infinitesimals were widely used before
Leibniz.  The novelty of Leibniz's law of continuity is not the
introduction of infinitesimals, but rather the formulation of a
resilient heuristic principle by Leibniz in an attempt to relate the
finite to the infinite.  Serfati's presentation of his epistemological
creationist principle as being on a higher level than Robinson's
transfer principle remains unsupported.  At any rate, what Serfati's
article \emph{does} reveal is that Serfati is perfectly aware of the
connection between Leibniz's law of continuity, on the one hand, and
the transfer principle for the hyperreals, on the other.

Two years ago, a new text by Serfati appeared as a book
chapter~\cite{Se10}.  Serfati reproduces an interesting Leibniz quote
on page~14 where Leibniz speaks explicitly of equality as encompassing
``equivalence up to an infinitely small'', the word ``equivalence"
originating in Leibniz.  Serfati proceeds to mention non-standard
analysis in the following terms:
\begin{quote}
contemporary mathematics has been able through non-standard analysis
to give a meaning (in a certain sense and at the cost of the
well-known complications) [to] the Leibnizian infinitely small
\cite[p.~15]{Se10}
\end{quote}
Serfati does {\em not\/} mention the fact that non-standard analysis
gives meaning to the law of continuity; he only mentions the fact that
it gives meaning to infinitesimals.

On page 26 Serfati discusses the disagreement between Cauchy and
Poncelet.  Poncelet exploited an extension of the principle of
continuity, known in the 18th and 19th centuries as the principle of
the ``generality of algebra", to study properties of algebraic curves,
and Cauchy was critical of this attitude, consistent with his
opposition to the ``generality of algebra'' in his {\em Cours
d'Analyse\/} and elsewhere (curiously, Serfati does not mention the
term ``generality of algebra").  Serfati proceeds to point out that
``Cauchy himself had failed in his work due to an intuitive
application--this time, false--of the schema" (of the principle of
continuity).  Here Serfati is referring to Cauchy's sum theorem of
1821, which he mentioned as an example (without using the name) a few
pages earlier.  Serfati is assuming that Cauchy made an error in 1821,
and ignores the essential ambiguity of the 1821 result (see Br\aa ting
\cite{Br}, Cutland et al. \cite{CKKR}).  Serfati is more explicit
about Cauchy's ``error" on page~28 and in footnote 33 (whose text
appears on page 30).  Serfati goes on to write:
\begin{quote}
Nowadays, therefore, 300 years after Leibniz, the ``principle of
continuity" belongs to an interiorized set of methodological rules.
Like the principle of symmetry (considered as normative) or that of
generalisation--exten\-sion (considered as a standard procedure of
construction of algebraic or topological objects), they constitute
part of the daily mathematical practice.  
\end{quote}
Serfati concludes:
\begin{quote}
Yet, they are never made explicit as such.  Somethimes there is talk
about this or that proof by continuity, but no manual discusses the
principle of such proofs \cite[p.~28]{Se10}.
\end{quote}

Now Keisler's \emph{instructor's manual} \cite{Ke76} (companion to his
{\em Elementary Calculus\/} \cite{Ke}) does discuss both Leibniz's law
of continuity, the transfer principle, and the relation between them
(see Subsection~\ref{42} above for some examples).  Thus, the
significance of Leibniz's heuristic principle, both philosophically
and mathematically, was and is clearly realized by the practitioners
of non-standard analysis.

Claiming that Leibniz's system for differential calculus was free of
logical fallacies may appear similar to claiming the possibility of
squaring the circle%
\footnote{Such was indeed the tenor of a recent referee report, see

\noindent
http://u.cs.biu.ac.il/$\sim$katzmik/straw2.html}
to a historian, but only if the latter internalized a triumviratist
spin on the history of mathematics as an ineluctable march, away from
logically fallacious infinitesimals, and toward the yawning heights of
Weierstrassian epsilontics.

\subsection{Cauchy and Moigno as seen by Schubring}
\label{cau}

To illustrate the significance of the distinction between Berkeley's
pair of criticisms, consider Cauchy's proof of the intermediate value
theorem \cite[Note~III, p.~460-462]{Ca21}.  Fowler~\cite{Fo} discussed
Cauchy's proof in the context of the decimal representation, and made
the connection between Cauchy and Stevin.%
\footnote{\label{IVT}Fowler notes that ``Stevin described an algorithm
for finding the decimal expansion of the root of any polynomial, the
same algorithm we find later in Cauchy's proof of the intermediate
value theorem'' \cite[p.~733]{Fo}.  The matter is discussed in detail
in B\l aszczyk et al. \cite{BKS}.  See also footnote~\ref{stevin}.}
Cauchy constructs an increasing sequence~$a_n$ and a decreasing
sequence~$b_n$ of successive approximations,~$a_n$ and~$b_n$ becoming
successively closer than any positive distance.%
\footnote{Cauchy's notation for the two sequences is~$x_0, x_1, x_2,
\ldots$ and~$X, X', X'', \ldots$ \cite[p.~462]{Ca21}.}
At this stage, the desired point is considered to have been exhibited,
by Cauchy.  A triumvirate scholar would object that Cauchy has not,
and could not have, proved the existence of the limit.  But imagine
that the {\em polytechnicien\/} Auguste Comte%
\footnote{Comte's notes of Cauchy's lectures have been preserved (see
\cite[p.~437]{Sch}).}
had asked {\em M.~le Professeur\/} Cauchy the following question:
\begin{quote}
Consider a decimal rank~$k>0$.  What is happening to the~$k$-th
decimal digit~$a_n^k$ of~$a_n$, and~$b_n^k$ of~$b_n$?
\end{quote}
Cauchy would have either sent Comte to the library to read Simon
Stevin%
\footnote{\label{stevin}Fearnley-Sander writes that ``the modern
concept of real number [...] was essentially achieved by Simon Stevin,
around 1600, and was thoroughly assimilated into mathematics in the
following two centuries'' \cite[p.~809]{Fea}.  Fearnley-Sander's
sentiment is echoed by van der Waerden \cite[p.~69]{van}.  Stevin had
anticipated Cauchy's proof of the intermediate value theorem, and
produced a fine-tuned version of the iteration, where each step of the
iteration produces an additional digit of the decimal expansion of the
solution.  The algorithm is discussed in more detail in \cite[\S10,
p.~475-476]{Ste}.  Stevin subdivides the interval into {\em ten\/}
equal parts, resulting in a gain of a new decimal digit of the
solution at every iteration of the algorithm.  Who needs the
``existence'' of the real numbers when Stevin constructs an explicit
decimal representation of the solution?  See also footnote~\ref{IVT}.}
(1548-1620), or else provided a brief argument to show that for~$n$
sufficiently large, the~$k$-th digit stabilizes, noting that special
care needs to be taken in the case when~$a_n$ is developing a tail
of~$9$s and~$b_n$ is developing a tail of~$0$s.%
\footnote{The pioneers of infinitesimal calculus were aware of the
non-uniqueness of decimal representation (at least) as early as 1770,
see Euler \cite[p.~170]{Eu}.}
Clearly the arguments appearing in Cauchy's book are sufficient to
identify the Stevin decimal expression of the limit.  From the modern
viewpoint, the only item missing is the remark that a Stevin decimal
{\em is\/} a number, {\em by definition\/} (modulo the identification
of the pair of tails).%
\footnote{Once the real numbers have been defined, considerable
technical difficulties remain in the definition of the multiplication
and other algebraic operations.  They were overcome by Dedekind (see
Fowler \cite{Fo}).}

Schubring reports on a purportedly successful effort by the cleric
Moigno, a student of Cauchy's, ``to pick apart'' (see
\cite[p.~445]{Sch}) infinitesimal methodologies.  Here Moigno puzzles
over how an infinitesimal magnitude can possibly be less than its own
half:
\begin{quote}
these magnitudes, [assumed to be] smaller than any given magnitude,
still have substance and are divisible[;~however,] their existence is
a chimera, since, necessarily greater than their half, their quarter,
etc., they are not actually less than any given magnitude (as quoted
by Schubring~\cite[p.~456]{Sch}).
\end{quote}
Moigno is clearly focusing on Berkeley's metaphysical criticism, not
the logical criticism.  As far as the metaphysical criticism is
concerned, the solution was readily available in Cauchy's work, and in
fact already in the work of Leibniz.  The solution is in terms of the
distinction between variable magnitude and constant magnitude; namely,
Cauchy's stratified hierarchy of constant quantities (A-continuum)
englobed inside an enriched B-continuum of variable quantities,%
\footnote{On occassion, Cauchy uses inequalities, rather than
equations involving infinitesimals, as, for instance, in Theorems I
and II in section 3 of chapter 2 of the \emph{Cours d'analyse} (see
Grabiner \cite{Gra}).  However, the thrust of his foundational
approach, {\em pace\/} Grabiner, is to use infinitesimals (generated
by null sequences) as inputs to functions, i.e., as individuals/atomic
entities; to define continuity in terms of infinitesimals; and to
apply infinitesimals to a range of problems, including an
infinitesimal definition of the ``Dirac'' delta function (see
Freudenthal \cite[p.~136]{Fr} and Laugwitz \cite{Lau92b}).}
including infinitesimals generated by null sequences (see
Section~\ref{cluster}).  

Schubring described Moigno as
\begin{quote}
the first writer to pick apart the traditional claim in favor of their
purported {\em simplicit\'e\/} \cite[p.~445]{Sch}.
\end{quote}
However, Moigno did not take apart the simplicity of infinitesimals,
either purported or real.  Rather, Moigno was confused, as many a
modern triumvirate scholar.  Similarly, in his 2007 anthology
\cite{Haw}, S.~Hawking reproduces Cauchy's {\em infinitesimal\/}
definition of continuity on page 639; but claims {\em on the same
page\/} 639, in a comic {\em non-sequitur\/}, that Cauchy ``was
particularly concerned to banish infinitesimals''.  In the same vein,
historian J.~Gray lists {\em continuity\/} among concepts Cauchy
allegedly defined
\begin{quote}
using careful, if not altogether unambiguous, {\bf limiting} arguments
\cite[p.~62]{Gray08} [emphasis added--authors],
\end{quote}
whereas in reality {\em limits\/} appear in Cauchy's definition only
in the sense of the {\em endpoints\/} of the domain of definition.

\subsection{Bishop and Connes}

Analyses of the critiques of Robinson's infinitesimals by E.~Bishop
and A.~Connes appear respectively in \cite{KK11d} and~\cite{KL}.

\section{\Horvath's analysis}
\label{horvath}

In a 1986 text in {\em Studia Leibnitiana\/}, M.~\Horvath~\cite{Ho86},
takes a critical view of a certain trend in contemporary Leibniz
scholarship.  We summarize some of \Horvath's main points below.

\subsection{Leibniz and Nieuwentijt}

Bernard Nieuwentijt possessed dramatically different intuitions about
infinitesimals as compared to Leibniz's.  Nieuwentijt favored
nilsquare infinitesimals: if~$dx$ is infinitesimal, then one should
have~$dx^2=0$.  Meanwhile, Leibniz was attached to
infinitesimals~$dx^n\not=0$ of arbitrary order.  With hindsight, we
know that both were right: Leibniz's intuitions are implemented in
Robinson's theory, whereas Nieuwentijt's, in Lawvere's Smooth
Infinitesimal Analysis (see J.~Bell \cite{Bel08, Bel09}).%
\footnote{\label{arthur2}It is therefore puzzling to find R.~Arthur
\cite{Ar} insisting, in J.~Bell's name, on similarities between
Leibniz's approach and that of Smooth Infinitesimal Analysis (SIA).
Arthur analyzes Bell's notion of an (intuitionistic) infinitesimal~$x$
as satisfying the relation~$\neg\neg\, x=0$.  Bell describes such an
$x$ as ``indistinguishable from 0''.  In more detail, Bell's
infinitesimal~$x$ satisfies NOT(NOT($x=0$)).  In classical logic this
would imply~$x=0$, but not in intuitionistic logic.  This can
certainly sound like a ``fictional" entity to a classically-trained
audience, but no intuitionist has been known to have embraced
fictionalism about anything in mathematics (certainly not
E.~Bishop--see \cite{KK11d}; indeed, many intuitionists insist that
all mathematical expressions refer to constructible objects), and at
any rate such ``fictionality" certainly has nothing to do with
Leibniz's.  Certainly Arthur's claim that SIA infinitesimals are
variable quantities unlike Robinson's is incorrect.  Arthur proceeds
to describe Bell's infinitesimals as ``fictional'', and bases his
analogy with Leibniz on the latter term, which procedure strikes us as
an unconvincing pun (see also footnote~\ref{arthur1}).}

What interests us here is Leibniz's response to Nieuwentijt's
criticism, discussed in \cite[p.~63]{Ho86}.  Nieuwentijt insists that
$dx^2=0$.  Leibniz responded in letters to both L'Hopital
\cite[p.~288]{Le95a} and to Huygens \cite[p.~207]{Le95b} dating from
1695, twenty years after the {\em Quadratura Arithmetica\/} (similar
comments in \emph{Cum Prodiisset} were discussed in
Subsection~\ref{41b}).

Leibniz responds that arbitrary orders of infinitesimals are necessary
so as to accomodate arbitrary orders of infinity (by inversion).  Note
that Nieuwentijt's criticism could have been answered more simply, by
affirming that~$dx$ is not an actual infinitesimal but merely a manner
of speaking, representing a shorthand for exhaustion \`a la
Archimedes.  Certainly Nieuwentijt did not hold that ordinary Stevin
numbers \cite{KK11c} can be nilsquare.  The fact that Leibniz does not
do so demonstrates that he and Nieuwentijt were of one mind as to the
non-Archimedean nature of~$dx$, and only disagreed about
``higher-order'' matters.

\subsection{Leibniz's 
pair of dual methodologies for justifying the calculus}

Leibniz pursued two different tracks for answering the critics of his
calculus: an Archimedean track and a track exploiting infinitesimals
with the law of continuity.  The latter track appears in texts
posterior to his {\em Quadratura Arithmetica\/}, e.g., texts dating
from 1680 and 1701 (see below). \Horvath{} writes:
\begin{quote}
The first approach is based on the fact that Leibniz regards his
infinitesimal calculus as an abbreviated language for proofs given by
the Greek method of exhaustion.  In connection with this method
Leibniz often refers to Archimedes by name
(\Horvath~\cite[p.~65-66]{Ho86}).
\end{quote}
We refer to this approach briefly as the A-methodology. \Horvath{}
proceeeds to describe the alternative methodology in the following
terms:
\begin{quote}
the second approach lies in the fact that, in Leibniz's mind, the
rules of his calculus can be proved by his so-called principle of
continuity \cite[p.~66]{Ho86}.
\end{quote}
The approach based on the law of continuity involves infinitesimals
(see Subsection~\ref{42} and Section~\ref{LOC}).  We refer to it
briefly as the B-methodology.

Leibniz explains the A-methodology in his letter to Pinson as follows:
\begin{quote}
In our calculations there is no need to conceive the infinite in a
rigorous way.  For instead of the infinite or the infinitely small,
one takes quantities as large, or as small, as necessary in order that
the error be smaller than the given error, so that one differs from
Archimedes' style only in the expressions, which are more direct in
our method and conform more to the art of invention (Leibniz
\cite{Le01a} cited in \Horvath~\cite[p.~66]{Ho86}).
\end{quote}
Leibniz clarifies his fictionalist position concerning the
B-methodology as involving ``ideal concepts'' in a 1702 letter to
Varignon in the following terms:
\begin{quote}
even if someone refuses to admit infinite and infinitely small lines
in a rigorous metaphysical sense and as real things, he can still use
them with confidence as ideal concepts which shorten the reasoning
(Leibniz \cite[p.~92]{Le02}, cited in \Horvath~\cite[p.~66]{Ho86}).
\end{quote}
Leibniz presented the dual pair of methodologies in an article dating
from the 1700s (see Subsection~\ref{105}).

\subsection{Incomparable quantities and the Archimedean property}

Leibniz uses the term ``incomparable quantities'' for quantities that
are supposed to implement the reduction of the calculus to the
A-methodology in a passage that also conveys the embattled situation
the defenders of the calculus found themselves in, in the face of
their ``opponents'':
\begin{quote}
These incomparable quantities are not at all fixed or determined but
can be taken to be as small as we wish in our geometrical reasoning
and so have the effect of the infinitely small in the rigorous sense.
If any opponent tries to contradict this proposition, it follows from
our calculus that the error will be less than any possible assignable
error, since it is in our power to take that incomparably small
quantity small enough for that purpose, inasmuch as we can always take
a quantity as small as we wish (Leibniz \cite[p.~92]{Le02} cited in
\Horvath~\cite[p.~66]{Ho86}).
\end{quote}
On the other hand, elsewhere Leibniz defines such ``incomparable
quantities'' in terms of the violation of what today is called the
Archimedean property.  Thus, Leibniz writes in a letter to
l'H\^opital:
\begin{quote}
I call incomparable quantities of which the one can not become larger
than the other if multiplied by any finite number.  This conception is
in accordance with the fifth definition%
\footnote{Leibniz is apparently referring to the \emph{fourth}
definition.}
of the fifth book of Euclid (Leibniz \cite[p.~288]{Le95a}, cited in
\Horvath~\cite[p.~63]{Ho86}).
\end{quote}
Here Leibniz employs the term ``incomparable quantity'' in the sense
of a non-Archimedean quantity, in the context of a B-continuum.

\subsection{From 
Leibniz's {\em Elementa\/} to his {\em Nova Methodus\/}}

Leibniz's manuscript {\em Elementa calculi novi\/} \cite{Le80} dates
from 1680.  The text {\em Elementa\/} is a preliminary draft of his
famous paper {\em Nova Methodus\/} \cite{Le84} dating from 1684.

The text had undergone a significant transformation between 1680 and
1684.  Thus, in {\em Elementa\/}, Leibniz describes quantities he
denotes 
\[
{}_1 D_2 C, {}_2C_3 D, \ldots
\]
as {\em ``incrementa momentanea''\/}, i.e., infinitely small.  By
1684, however, a change occurs:

\begin{quote}
Leibniz does not use differentials but only differences in the sense
of fixed, small, finite quantities.  Leibniz [presumably] does not use
the term ``infinitely small'' in his article [1684] in order to avoid
controversies which most likely would have arisen in connection with
this notion (\Horvath~\cite[p.~62]{Ho86}).
\end{quote}
Thus, in 1675 and 1684, Leibniz emphasizes the A-methodology, but five
years after {\em Quadratura Arithmetica\/} and four years before {\em
Nova Mathodus\/}, he relies upon the B-methodology.

\subsection{{\em Justification du calcul des infinitesimales\/}}
\label{105}

Leibniz penned an additional defense of his calculus in the early
1700s.  The 1700 version is entitled {\em Defense du calcul des
differences\/} \cite{Le00}.  The 1701 version is entitled {\em
Justification du Calcul des infinitesimales\ldots}~\cite{Le01b}, and
is the published version of the article.  In his {\em Defense\/},
\begin{quote}
Leibniz appeals to the method of Archimedes in the first part of the
draft, while the allusion to the law of continuity can be found in the
second part of his draft \cite[footnote~27, p.~69]{Ho86}.
\end{quote}
The order is reversed in his {\em Justification\/}.  Thus, both the
A-methodology and the B-methodology are present at this late stage,
and it is the latter that gets right-of-way, a quarter century after
{\em Quadratura Arithmetica\/}.

\section{A Leibniz--Robinson route out of the labyrinth?}
\label{conclusion}

We have dissected Berkeley's critique into its component parts
following Sherry, and have revealed the implausibility of some of the
assumptions underlying that critique.  We have discussed both the
critique's ill-informed nature, and Berkeley's contradictory attitude
when writing about a different field of mathematics, such as
arithmetic.

A significant, and widely denied, aspect of the story is the existence
of a direct perceptual, cognitive, and even formal connection between
historical infinitesimals as they were practiced by giants like
Leibniz and Cauchy, via the work of Stolz, Paul du Bois-Reymond,
Veronese, and others at the turn of the century, and the emergence of
hyper-real fields in the middle of the 20th century.

Berkeley's {\em Philosophical Commentaries\/} \cite{Be48} is an early
work, dating from approximately 1709.  Here Berkeley asserts that
Euclidean geometry is full of paradoxes.  He backs away from this
assertion as he matures, but later in his \emph{Principles}, his most
mature work, it turns out that he won't accept infinite divisibility
and sees it on a par with infinitesimals.  How is it that Berkeley has
no difficulty with imaginary roots and yet balks at infinite
divisibility?

Berkeley is prepared to be instrumentalist about imaginary roots, yet
refuses to be instrumentalist about certain idealisations of the
continuum involved in the calculus.  Another puzzling aspect of his
position is the following.  If he indeed does not accept the latter,
why does he provide a ``cancellation of errors" justification of the
latter?  The justification turns out to be meaningless symbol-pushing,
according to a recent article by Andersen \cite{An11}; or more
precisely \emph{circular} symbol pushing: the fact that everything
works out in the end is based on a result due to Apollonius of Perga
\cite[Book I, Theorem 33]{Ap} on the tangents of the parabola, which
is equivalent to the calculation of the derivative in the quadratic
case.  Thus, Berkeley's logical criticism is further weakened by the
circular \emph{logic} he relied upon in his ``compensation of errors''
approach.  Given Berkeley's fame among historians of mathematics for
allegedly spotting logical flaws in infinitesimal calculus, it is
startling to spot circular logic at the root of Berkeley's own
doctrine of compensation of errors, seeing that his new, improved
calculation of the derivative of~$x^2$ relies upon Apollonius'
determination of the tangent to a parabola.

There is indeed an irreducible conceptual clash between Leibniz and
Berkeley.  One of the sides had to give way.  What historians
sometimes do not fully appreciate is the fact that the weaker side was
Berkeley's, not Leibniz's.  Modern mathematics would not even start
without non-referential concepts that would have been ridiculed by
Berkeley as meaningless due to his empiricist bias.

Edwin Hewitt \cite{Hew} in 1948 constructed infinitesimal-enriched
continua and introduced the term {\em hyper-real\/}.  In 1955,
J.~\Los~\cite{Lo55} proved what has come to be known as \Los's theorem
(for ultraproducts), whose consequence is the transfer principle for
hyper-real fields, which embodies a mathematical implementation of
Leibniz's heuristic law of continuity.

\section*{Acknowledgements}

We are grateful to H.~Jerome Keisler for helpful remarks that helped
improve an earlier version of the manuscript.  The influence of Hilton
Kramer (1928-2012) is obvious.

\appendix

\section{Rival continua}
\label{rival}

The historical roots of infinitesimals go back to Cauchy, Leibniz, and
ultimately to Archimedes.  Cauchy's approach to infinitesimals is not
a variant of the hyperreals.  Rather, Cauchy's work on the rates of
growth of functions anticipates the work of late 19th century
investigators such as Stolz, du Bois-Reymond, Veronese, Levi-Civita,
Dehn, and others, who developed non-Archimedean number systems against
virulent opposition from Cantor, Russell, and others (see Ehrlich
\cite{Eh06} and Katz and Katz \cite{KK11a} for details).  The work on
non-Archimedean systems motivated the work of T.~Skolem on
non-standard models of arithmetic \cite{Sk}, which subsequently
stimulated work culminating in the hyperreals of Hewitt, \L os, and
Robinson.

\begin{figure}
\includegraphics[height=2in]{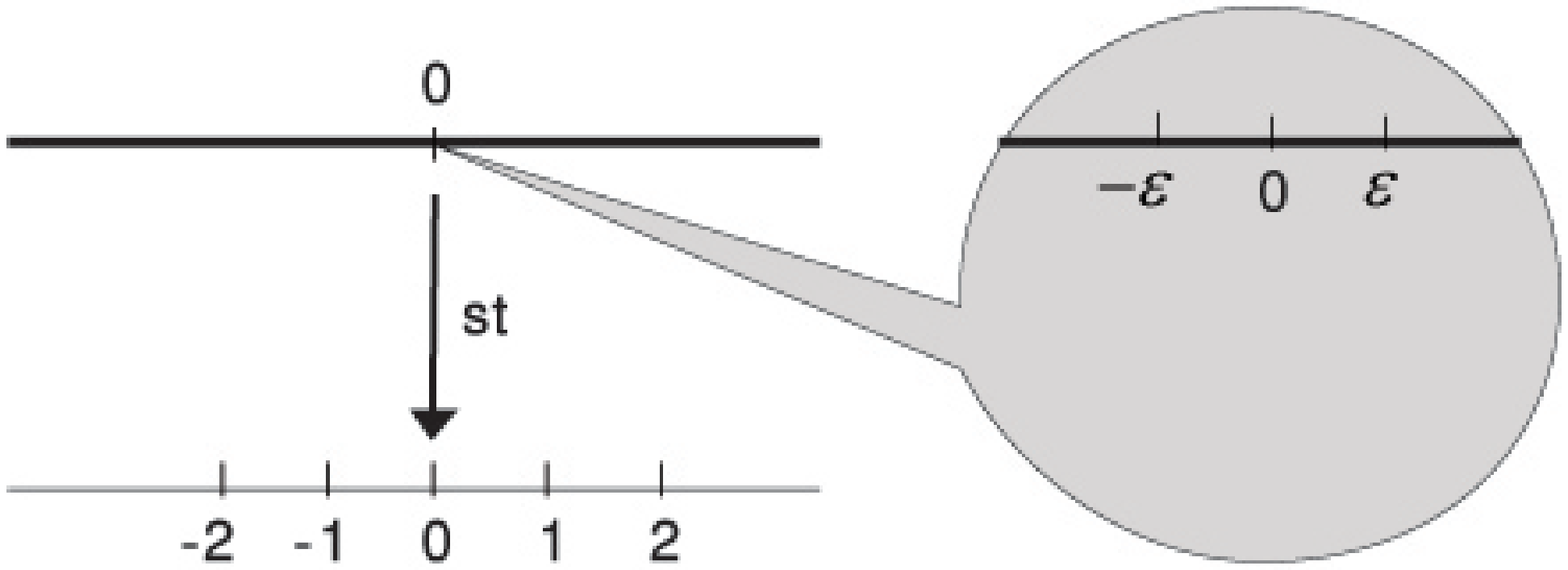}
\caption{\textsf{Zooming in on infinitesimal~$\epsilon$ (here st$(\pm
\epsilon)=0$)}}
\label{tamar}
\end{figure}

The relation between the rival theories of the continuum distinguished
by Felix Klein (see Subsection~\ref{klein}) can be summarized as
follows.  A Leibnizian definition of the differential quotient
\[
\frac{\Delta y}{\Delta x},
\] 
whose logical weakness was criticized by Berkeley, was modified by
A.~Robinson by exploiting a map called {\em the standard part\/},
denoted~``st'', from the finite part of a ``thick'' B-continuum (i.e.,
a Bernoullian continuum),%
\footnote{\label{f61}See footnote~\ref{f49}.}
to a ``thin'' A-continuum (i.e., an Archimedean continuum), as
illustrated in Figures~\ref{31} and \ref{tamar}.  The derivative is
defined as~$\hbox{st} \left( \frac{\Delta y}{\Delta x} \right)$,
rather than the differential quotient ~$\frac{\Delta y}{\Delta x}$
itself.  Robinson wrote that ``this is a small price to pay for the
removal of an inconsistency'' (Robinson \cite[p.~266]{Ro66}).
However, the process of discarding the higher-order infinitesimals has
solid roots in Leibniz's law of homogeneity (see
Subsection~\ref{homogeneity}).

\subsection{Hyperreals via maximal ideals}

We summarize a 20th century implementation of an alternative to an
Archimedean continuum, namely an infinitesimal-enriched continuum.
Such a continuum is not to be confused with incipient notions of such
a continuum found in earlier centuries.  We refer to such a continuum
as a B-continuum.$^{\ref{f61}}$.  We begin with a heuristic
representation of a B- or ``thick'' continuum, denoted~$\RRR$, in
terms of an infinite resolution microscope (see Figure~\ref{tamar}).
One presentation of such a structure is in Robinson \cite{Ro61}.  Such
an infinitesimal-enriched continuum is suitable for use in calculus,
analysis, and elsewhere.  Robinson built upon earlier work by
E.~Hewitt \cite{Hew}, J.~\Los{} \cite{Lo55}, and others.  In 1962,
W.~Luxemburg \cite{Lu62} popularized a presentation of Robinson's
theory in terms of the ultrapower construction,%
\footnote{Note that both the term ``hyper-real'', and an ultrapower
construction of a hyperreal field, are due to E.~Hewitt in 1948, see
\cite[p.~74]{Hew}.  Luxemburg~\cite{Lu62} also clarified its relation
to the competing construction of Schmieden and Laugwitz \cite{SL},
also based on sequences, which used a different kind of filter.}
in the mainstream foundational framework of the Zermelo--Fraenkel set
theory with the axiom of choice (ZFC).

The construction can be viewed as a relaxing, or refining, of Cantor's
construction of the reals.  This can be motivated by a discussion of
rates of convergence as follows.  In Cantor's construction, a real
number $u$ is represented by a Cauchy sequence $\langle u_n :
n\in\N\rangle$ of rationals.  But the passage from $\langle
u_n\rangle$ to $u$ in Cantor's construction sacrifices too much
information.  We would like to retain a bit of the information about
the sequence, such as its ``speed of convergence".  This is what one
means by ``relaxing" or ``refining" Cantor's construction of the reals
(cf.~Giordano et al.~\cite{GK11}).  When such an additional piece of
information is retained, two different sequences, say $\langle
u_n\rangle$ and $\langle u'_n\rangle$, may both converge to $u$, but
at different speeds.  The corresponding ``numbers" will differ from
$u$ by distinct infinitesimals.  If $\langle u_n\rangle$ converges to
$u$ faster than $\langle u'_n\rangle$, then the corresponding
infinitesimal will be smaller.  The retaining of such additional
information allows one to distinguish between the equivalence class of
$\langle u_n\rangle$ and that of $\langle u'_n\rangle$ and therefore
obtain distinct hyperreals infinitely close to $u$.

At the formal level, we proceed as follows.  We construct a hyperreal
field as a quotient of the collection of arbitrary sequences, where a
sequence
\begin{equation}
\langle u_1, u_2, u_3, \ldots \rangle
\end{equation}
converging to zero generates an infinitesimal (the kernel of the
quotient homomorphism is the maximal ideal $\mathcal{M}$ described
below).  Arithmetic operations are defined at the level of
representing sequences; e.g., addition and multiplication are defined
term-by-term.  Thus, we start with the ring~$\Q^\N$ of sequences of
rational numbers.  Let
\begin{equation}
\label{22new}
\mathcal{C}_\Q \subset \Q^\N
\end{equation}
denote the subring consisting of Cauchy sequences.  The reals are by
definition the quotient field
\begin{equation}
\label{realbis}
\R:= \mathcal{C}_\Q / \fnull,
\end{equation}
where~$\fnull$ is the ideal containing all null sequences (i.e.,
sequences tending to zero).%
\footnote{Namely, the traditional construction of the real field,
usually attributed to Cantor, views a real number as an equivalence
class of Cauchy sequences of rational numbers.  Null sequences
comprise the equivalence class corresponding to the real number
$0\in\R$.}
Note that $\Q$ is imbedded in $\Q^\N$ by constant sequences.  An
infinitesimal-enriched extension of~$\Q$ may be obtained by
modifying~\eqref{realbis}.  Now consider the subring
\[
\fez\subset\fnull
\]
of sequences that are ``eventually zero'', i.e., vanish at all but
finitely many places.  Then the quotient~$\mathcal{C}_\Q / \fez$
naturally surjects onto~$\R= \mathcal{C}_\Q / \fnull$.  The elements
in the \emph{kernel} of the surjection
\[
\mathcal{C}_\Q/\fez \to \R
\]
are prototypes of infinitesimals.%
\footnote{Such elements could be called ``infinitesimal'' to the
extent that they violate the Archimedean property suitably
interpreted, but the ring they are elements of has unsatisfactory
properties.}
Note that the quotient~$\mathcal{C}_\Q/\fez$ is not a field, as~$\fez$
is not a maximal ideal.  To obtain a field, we must replace~$\fez$ by
a \emph{maximal} ideal.

It is more convenient to describe the modified construction using the
ring~$\R^\N$ rather than~$\mathcal{C}_\Q$ of~\eqref{22new}.

We therefore redefine~$\fez$ to be the ring of \emph{real} sequences
in~$\R^\N$ that eventually vanish, and choose a \emph{maximal} proper
ideal~$\mathcal{M}$ so that we have
\begin{equation}
\label{23new}
\fez\subset\mathcal{M}\subset\R^\N.
\end{equation}
Then the quotient
\begin{equation}
\label{A5}
\RRR:=\R^\N/\mathcal{M}
\end{equation}
is a hyperreal field.  The foundational material needed to ensure the
existence of a maximal ideal~$\mathcal{M}$ satisfing~\eqref{23new} is
weaker than the axiom of choice.  This concludes the construction of a
hyperreal field~$\RRR$ in the traditional foundational framework, ZFC.

The construction is equivalent to the usual ultrapower construction as
popularized by Luxemburg.%
\footnote{\label{f65}Analyzing $\mathcal{M}$, one discovers that its
structure is controlled by a free ultrafilter on $\N$.  The order
relation on $\RRR$ is defined relative to the ultrafilter.  Additional
details may be found in \cite[Appendix~A]{BKS}.}
Thus it is not entirely accurate to suppose, as Jesseph does, that a
consistent theory of infinitesimals requires the resources of model
theory.  The resources of a rigorous undergraduate course in abstract
algebra suffice.

\subsection{Example}

To give an example, the sequence
\begin{equation}
\label{infinitesimal}
\left\langle \tfrac{1}{n} : n\in \N\right\rangle
\end{equation}
represents a nonzero infinitesimal, in the sense that its class
$\left[\tfrac{1}{n} \right]$ in~\eqref{A5} is nonzero and satisfies
$\left[\tfrac{1}{n} \right]<r$ for every positive real number $r$.%
\footnote{See footnote~\ref{f65}.}

\subsection{Construction of standard part}
\label{A2}

In the field~$\RRR$ of~\eqref{A5}, consider the subring~$I\subset\RRR$
consisting of infinitesimal elements (i.e., elements~$e$ such
that~$|e|<\frac{1}{n}$ for all~$n\in\N$).  Denote by~$I^{-1}$ the set
of inverses of nonzero elements of~$I$.  The complement~$\RRR \,
\setminus \, I^{-1}$ consists of all the finite (sometimes called
\emph{limited}) hyperreals.  Constant sequences provide an
inclusion~$\R\subset\RRR$.  Every element~$x\in \RRR\setminus I^{-1}$
is infinitely close to some real number~$x_0\in\R$.  The
\emph{standard part function}, denoted ``st'', associates to every
finite hyperreal, the unique real infinitely close to it:
\[
\st:\RRR\setminus I^{-1} \to \R, \text{\;with\;\;} x \mapsto x_0.
\]
The real~$x_0$ is sometimes called the \emph{shadow} of~$x$.  If $x$
happens to be the equivalence class of a Cauchy sequence $\langle
x_n:n\in\N \rangle$, then the shadow of $x$ is the limit of $\langle
x_n\rangle$:
\[
\st(x)=\lim_{n\to\infty} x_n.
\]
As explained in Subsection~\ref{homogeneity}, the standard part
function can be seen as an implementation of Leibniz's transcendental
law of homogeneity.

\subsection{The transfer principle}

The {\em transfer principle\/} is a mathematical implementation of
Leibniz's heuristic {\em law of continuity\/}: ``what succeeds for the
finite numbers succeeds also for the infinite numbers and vice versa''
(see Robinson~\cite[p.~266]{Ro66}).  The transfer principle, allowing
an extention of every first-order real statement%
\footnote{Some examples were provided in Subsection~\ref{45} following
formula~\eqref{4.4}.}
to the hyperreals, is a consequence of the theorem of J.~{\L}o{\'s} in
1955 (see~\cite{Lo55}), and can therefore be referred to as a
Leibniz-{\L}o{\'s} transfer principle.  A Hewitt-{\L}o{\'s} framework
allows one to work in a B-continuum satisfying the transfer principle.

\begin{figure}
\includegraphics[height=1.6in]{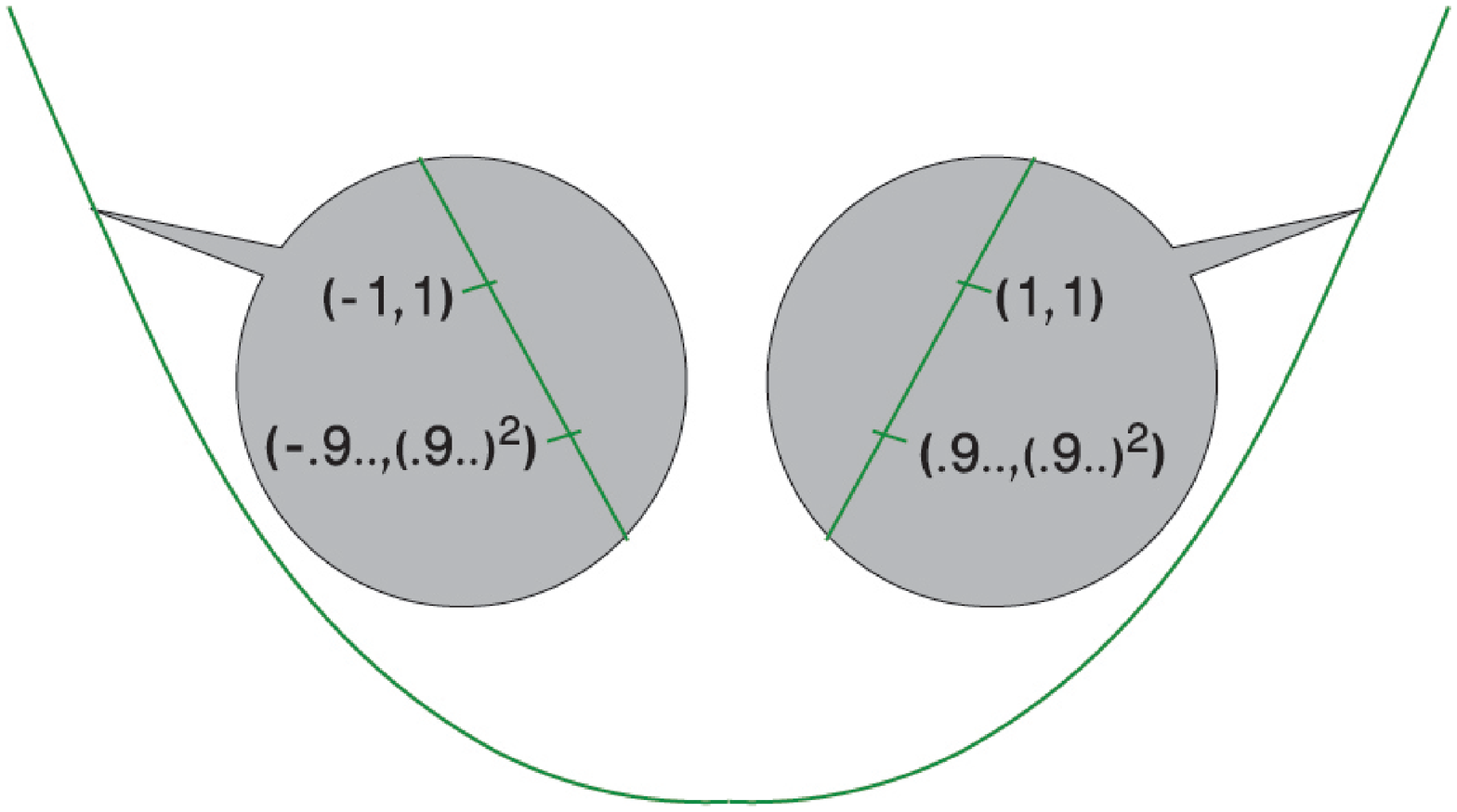}
\caption{\textsf{Differentiating~$y=f(x)=x^2$ at~$x=1$ yields
$\tfrac{\Delta y}{\Delta x} = \tfrac{f(.9..) - f(1)}{.9..-1} =
\tfrac{(.9..)^2 - 1}{.9..-1} = \tfrac{(.9.. - 1)(.9.. + 1)}{.9..-1} =
.9.. + 1 \approx 2.$ Here~$\approx$ is the relation of being
infinitely close.  Hyperreals of the form~$.9..$ are discussed in
\cite{KK2}}}
\label{jul10}
\end{figure}

A helpful ``semicolon'' notation for presenting an extended decimal
expansion of a hyperreal was described by A.~H.~Lightstone~\cite{Li}.
See also P.~Roquette \cite{Roq} for infinitesimal reminiscences.  A
discussion of infinitesimal optics is in K.~Stroyan \cite{Str},
J.~Keisler~\cite{Ke}, D.~Tall~\cite{Ta80}, L.~Magnani \&
R.~Dossena~\cite{MD, DM}, and Bair \& Henry~\cite{BH}.  Applications
of the B-continuum range from aid in teaching calculus \cite{El, KK1,
KK2, Ta91, Ta09a} (see illustration in Figure~\ref{jul10}) to the
Bolzmann equation (see L.~Arkeryd~\cite{Ar81, Ar05}); modeling of
timed systems in computer science (see H.~Rust \cite{Rust}); Brownian
motion and economics (see Anderson \cite{An76}); mathematical physics
(see Albeverio {\em et al.\/} \cite{Alb}); etc.  The hyperreals can be
constructed out of integers (see Borovik, Jin, \& Katz \cite{BJK}).
The traditional quotient construction using Cauchy sequences, usually
attributed to Cantor, can be factored through the hyperreals (see
Giordano \& Katz \cite{GK11}).

\bigskip

\textbf{Mikhail G. Katz} is Professor of Mathematics at Bar Ilan
University, Ramat Gan, Israel.  His joint study with P.~B\l aszczyk
and D.~Sherry entitled \emph{Ten misconceptions from the history of
analysis and their debunking} is due to appear in \emph{Foundations of
Science}.  A joint study with A.~Borovik entitled \emph{Who gave you
the Cauchy--Weierstras tale?  The dual history of rigorous calculus}
appeared in \emph{Foundations of Science} (online first).  A joint
study with A.~Borovik and R.~Jin entitled \emph{An integer
construction of infinitesimals: Toward a theory of Eudoxus hyperreals}
is due to appear in \emph{Notre Dame Journal of Formal Logic}
\textbf{53} (2012), no.~4.  A joint study with Karin Katz entitled
\emph{A Burgessian critique of nominalistic tendencies in contemporary
mathematics and its historiography}, appeared in \emph{Foundations of
Science} \textbf{17} (2012), no.~1, 51--89.  A joint study with Karin
Katz entitled \emph{Meaning in classical mathematics: is it at odds
with Intuitionism?} appeared in {\em Intellectica\/} \textbf{56}
(2011), no.~2, 223-302.  A joint study with Karin Katz, \emph{Stevin
numbers and reality}, appeared in \emph{Foundations of Science}
(online first).  A joint study with Eric Leichtnam, \emph{Commuting
and non-commuting infinitesimals}, is due to appear in \emph{American
Mathematical Monthly}.  A joint study with David Tall entitled
\emph{The tension between intuitive infinitesimals and formal
mathematical analysis}, appeared as a chapter in a book edited by
Bharath Sriraman, \emph{Crossroads in the History of Mathematics and
Mathematics Education}, {\em The Montana Mathematics Enthusiast
Monographs in Mathematics Education\/} \textbf{12}, Information Age
Publishing, Inc., Charlotte, NC, 2012.

\bigskip

{\bf David Sherry} is fortunate to be professor of philosophy at
Northern Arizona University in the cool pines of the Colorado Plateau.

\end{document}